\documentclass[a4paper, 11pt]{amsart}

\usepackage{cite}
\usepackage{hyperref}
\usepackage{graphicx}
\usepackage{color}
\usepackage[T1]{fontenc}
\usepackage{amsmath}
\usepackage{amssymb}
\usepackage{amsthm}
\usepackage{bm}
\usepackage{eucal}
\usepackage{braket}
\usepackage{stmaryrd}
\usepackage{mathrsfs}
\usepackage[all,2cell]{xy}
\UseAllTwocells
\usepackage{youngtab}
\usepackage{tikz}
\usetikzlibrary{intersections,calc,arrows.meta}
\calclayout

\newcommand{\sgn}{\mathrm{sgn}}

\newcommand{\pr}[1]{#1^{\prime}}
\newcommand{\ppr}[1]{#1^{\prime\prime}}

\newcommand{\what}[1]{\widehat{#1}}
\newcommand{\wtilde}[1]{\widetilde{#1}}

\newcommand{\End}{\mathrm{End}}

\newcommand{\Hom}{\mathrm{Hom}}

\newcommand{\tr}{{\mathrm{tr}}}

\newcommand{\id}{\mathrm{id}}

\newcommand{\catMod}{\mathsf{Mod}}

\newcommand{\funInd}{\mathsf{Ind}}
\newcommand{\funRes}{\mathsf{Res}}
\newcommand{\funId}{\mathsf{Id}}

\newcommand{\hackcenter}[1]{
 \xy (0,0)*{#1}; \endxy}

\newcommand{\bC}{\mathbb{C}}

\newcommand{\bR}{\mathbb{R}}

\newcommand{\bY}{\mathbb{Y}}
\newcommand{\bZ}{\mathbb{Z}}

\newcommand{\bfB}{\mathbf{B}}
\newcommand{\bfm}{\mathbf{m}}

\newcommand{\cC}{\mathcal{C}}

\newcommand{\cE}{\mathcal{E}}
\newcommand{\cF}{\mathcal{F}}

\newcommand{\cH}{\mathcal{H}}

\newcommand{\cP}{\mathcal{P}}

\newcommand{\sfm}{\mathsf{m}}

\makeatletter
\newcommand{\xRightarrow}[2][]{\ext@arrow 0359\Rightarrowfill@{#1}{#2}}
\makeatother

\theoremstyle{plain}
\newtheorem{thm}{Theorem}[section]
\newtheorem{lem}[thm]{Lemma}
\newtheorem{prop}[thm]{Proposition}
\newtheorem{cor}[thm]{Corollary}
\theoremstyle{definition}
\theoremstyle{remark}
\newtheorem{rem}[thm]{Remark}%[section]

\makeatletter
    
    \@addtoreset{equation}{section}
\makeatother

\title[Normalized characters and Boolean cumulants]{Normalized characters of symmetric groups and Boolean cumulants via Khovanov's Heisenberg category}
\author{Shinji Koshida}
\address{Department of Mathematics and Systems Analysis, Aalto University}
\email{shinji.koshida@aalto.fi}

\begin{document}

\begin{abstract}
%The normalized characters of symmetric groups are known to be well-behaved when their asymptotics are concerned.
In this paper, we study relationships between the normalized characters of symmetric groups and the Boolean cumulants of Young diagrams.
Specifically, we show that each normalized character is a polynomial of twisted Boolean cumulants with coefficients being non-negative integers,
and conversely, that, when we expand a Boolean cumulant in terms of normalized characters, the coefficients are again non-negative integers.
The main tool is Khovanov's Heisenberg category and the recently established connection of its center to the ring of functions on Young diagrams,
which enables one to apply graphical manipulations to the computation of functions on Young diagrams.
Therefore, this paper is an attempt to deepen the connection between the asymptotic representation theory and graphical categorification.
\end{abstract}

\maketitle

%\tableofcontents

\section{Introduction}
\subsection{Backgrounds}
Study of the asymptotic behavior of representations of symmetric groups as the group size goes to infinity was initiated by Vershik and Kerov~\cite{VershikKerov1981}, and has been developed as a central subject of the so-called asymptotic representation theory.
The asymptotic representation theory has not only become a field that gives general ideas of studying towers of algebraic structures such as compact groups~\cite{Voiculescu1976,VershikKerov1982,Biane1995,CollinsSniady2009,BufetovGorin2015,CollinsNovakSniady2018}, Hecke algebras~\cite{FerayMeliot2012}, quantum groups~\cite{Sato2019,Sato2021a,Sato2021b,Sato2021c}, among others apart from symmetric groups~\cite{Biane1998,Kerov2003,BorodinOlshanski_book2017},
but also caused a new stream of probability theory, integrable probability (for survey, see e.g. \cite{BorodinGorin2012,BorodinPetrov2014,Corwin2014}),
which has now been widening its perspective even independently of representation theory.

% Apart from symmetric groups, several different algebraic structures including unitary groups, Hecke algebras, quantum groups, are studied in asymptotic representation.
% Furthermore asymptotic representation had an impact on integrable probability.

In the case of symmetric groups, the irreducible representations are labelled by Young diagrams (see e.g.~\cite{Fulton1996,Sagan2001}).
One approach to the asymptotic representation of symmetric groups is
to study the asymptotic behavior of characters viewed as functions of Young diagrams.
As a result, several observables defined on Young diagrams show up to describe the asymptotic behavior.
Most prominently, Biane~\cite{Biane1998} pointed out that free cumulants play significant roles in the asymptotic representation of symmetric groups and opened up a new approach from free probability and random matrix theory~\cite{HiaiPetz2000,NicaSpeicher2006,MingoSpeicher2017}.
In this paper, however, we are concerned with another series of observables, Boolean cumulants~\cite{SpeicherWoroudi1993}.

Another pillar of this paper is Khovanov's Heisenberg category introduced in \cite{Khovanov2014} for the purpose of categorifying the infinite dimensional Heisenberg algebra.
% Khovanov's Heisenberg category and its quantum deformation have already had many applications in categorification of algebraic structures and algebraic results.
It is a category whose morphisms are given by planar diagrams modulo local relations that reflect the induction and restriction of representations of symmetric groups. It was already anticipated in \cite{Khovanov2014} that Khovanov's Heisenberg category could be related to free probability theory in the spirit of \cite{GuionnetJonesShlyakhtenko2010}, which emphasized that the planar structure is shared by subfactors from operator algebras, free probability and random matrices.
Remarkably, in this direction, \cite{KLM2019} established an isomorphism between the center of Khovanov's Heisenberg category and the ring of shifted symmetric functions,
and in particular clarified that elements of the center of Khovanov's Heisenberg category give rise to observables on Young diagrams.

We view the results of \cite{KLM2019} as what allows for graphical calculus of observables on Young diagrams.
In the rest of the paper, we investigate the relation between normalized characters of symmetric groups and Boolean cumulants
standing on this viewpoint.

% The main goal in this text is to prove \cite[Conjecture 8]{RS2008}; that is
% \begin{thm}
% For a normalized character $\Sigma_{k_{1},\dots,k_{l}}$,
% its sign change $(-1)^{l}\Sigma_{k_{1},\dots,k_{l}}$ is a polynomial of the twisted Boolean cumulants $\hat{B}_{i}$ with coefficients being non-negative integers.
% \end{thm}

% The key tool was presented in \cite{KLM2019}, where the center of Khovanov's Heisenberg category and the ring of shifted symmetric functions were identified.
% This identification enables us to carry out diagrammatic calculations of normalized characters and Boolean cumulants.

% \section{Preliminaries and main result}
\subsection{Symmetric groups}
For $n\in \bZ_{\geq 0}$, let $S_{n}$ be the $n$-th symmetric group (under the convention that $S_{0}=\{e\}$ is the trivial group).
We write $\bY_{n}$ for the set of Young diagrams with $n$ cells (under the convention that $\bY_{0}=\{\emptyset\}$).
Then the isomorphism classes of irreducible representations of $S_{n}$ are labelled by $\bY_{n}$;
for $\lambda\in \bY_{n}$, let $(\rho^{\lambda},V^{\lambda})$ be an irreducible representation corresponding to $\lambda$.
See e.g.~\cite{Fulton1996} for accounts of representations of symmetric groups (or~\cite{OkounkovVershik1996,VershikOkounkov2005} for an alternative approach.)
We define the character of this representation by
\begin{equation*}
    \chi^{\lambda}(\sigma)=\frac{\tr_{V^{\lambda}}\rho^{\lambda}(\sigma)}{\tr_{V^{\lambda}}\rho^{\lambda}(e)},\quad \sigma\in S_{n},
\end{equation*}
where $e\in S_{n}$ is the unit element.
Note that in this convention, we have $\chi^{\lambda}(e)=1$.
The definition is indeed independent of the choice of an irreducible representation of type $\lambda$.
The character $\chi^{\lambda}$ only depends on conjugacy classes of $S_{n}$,
which are in one-to-one correspondence with the set $\cP_{n}$ of partitions of $n$.
Although the set $\cP_{n}$ of partitions is naturally identified with the set of Young diagrams $\bY_{n}$, we distinguish them because of their different roles.
We write $\chi^{\lambda}_{\pi}$ for the value of $\chi^{\lambda}$ on the conjugacy class corresponding to a partition $\pi\in \cP_{n}$.

For a fixed partition $\pi\in \cP_{k}$, the normalized character $\Sigma_{\pi}$ is the function on $\bY:=\bigcup_{n=0}^{\infty}\bY_{n}$ defined by
\begin{equation*}
    \Sigma_{\pi}(\lambda):=
    \begin{cases}
        (n\downharpoonright k)\chi^{\lambda}_{\pi\cup (1^{n-k})},& \lambda\in \bY_{n},\, n\geq k, \\
        0 & \mbox{otherwise} .
    \end{cases}
\end{equation*}
Here $(n\downharpoonright k)=n(n-1)\cdots (n-k+1)$ is the falling factorial,
and $\pi\cup (1^{n-k})$ is the partition of $n$ obtained by adding $n-k$ copies of $1$ to $\pi$.

Collecting all normalized characters $\Sigma_{\pi}$ for $\pi\in \cP:=\bigcup_{n=0}^{\infty}\cP_{n}$ with the convention $\cP_{0}=\{\emptyset\}$,
we obtain a family of functions $\{\Sigma_{\pi}:\pi\in \cP\}$ on $\bY$.

\subsection{Boolean cumulants}
\begin{figure}
\begin{center}
\begin{tikzpicture}
    \draw (-6.5,0)--(6.5,0);
    \draw (-6.5,6.5)--(0,0)--(6.5,6.5);
    \draw[very thick] (-6.5,6.5)--(-5,5)--(-4,6)--(-3,5)--(-2,6)--(0,4)--(1,5)--(2,4)--(4,6)--(5,5)--(6.5,6.5);
    \draw (-4,4)--(-3,5);
    \draw (-3,3)--(-1,5);
    \draw (-2,2)--(0,4);
    \draw (-1,1)--(2,4);
    \draw (1,1)--(-3,5);
    \draw (2,2)--(0,4);
    \draw (3,3)--(2,4);
    \draw (4,4)--(3,5);
    \draw[dashed] (-5,5)--(-5,0)node[below]{$x_{1}$};
    \draw[dashed] (-3,5)--(-3,0)node[below]{$x_{2}$};
    \draw[dashed] (0,4)--(0,0)node[below]{$x_{3}$};
    \draw[dashed] (2,4)--(2,0)node[below]{$x_{4}$};
    \draw[dashed] (5,5)--(5,0)node[below]{$x_{5}$};
    \draw[dashed] (-4,6)--(-4,0)node[below]{$y_{1}$};
    \draw[dashed] (-2,6)--(-2,0)node[below]{$y_{2}$};
    \draw[dashed] (1,5)--(1,0)node[below]{$y_{3}$};
    \draw[dashed] (4,6)--(4,0)node[below]{$y_{4}$};
\end{tikzpicture}
\caption{The profile of $\lambda=(5,3,2,2,1)$. The local minima and maxima are given by $(x_{1},x_{2},x_{3},x_{4},x_{5})=(-5,-3,0,2,5)$ and $(y_{1},y_{2},y_{3},y_{4})=(-4,-2,1,4)$.}
\label{fig:profile_example}
\end{center}
\end{figure}
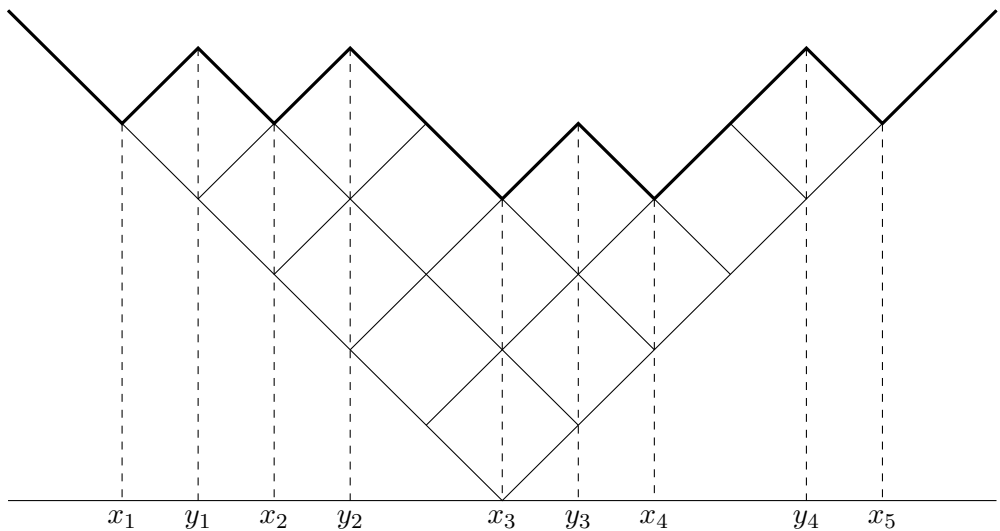

We employ the Russian convention to draw Young diagrams.
In this convention, a Young diagram is identified with a function $\omega$ on $\bR$ such that
\begin{enumerate}
    \item $\omega (x)=|x|$ when $|x|$ is sufficiently large,
    \item $\omega$ is piecewise linear with slopes being either $1$ or $-1$, and
    \item local maxima and minima are located on $\bZ$.
\end{enumerate}
Such a function is referred to as the profile of a Young diagram,
and it can be recovered from its local minima $x_{1}<x_{2}<\cdots<x_{n}$
and local maxima $y_{1}<y_{2}<\cdots <y_{n-1}$.
Note that the local minima and maxima of a profile are always interlacing:
\begin{equation}
\label{eq:interlacing_xy}
    x_{1}<y_{1}<x_{2}<y_{2}<\cdots<y_{n-1}<x_{n},
\end{equation}
and that their sums coincide:
\begin{equation}
\label{eq:sums_xy}
    \sum_{i=1}^{n}x_{i}=\sum_{i=1}^{n-1}y_{i}.
\end{equation}
For instance, the profile of $\lambda=(5,3,2,2,1)$ is depicted in Figure~\ref{fig:profile_example}.
In this example, the local minima and maxima are given by $(x_{1},x_{2},x_{3},x_{4},x_{5})=(-5,-3,0,2,5)$ and $(y_{1},y_{2},y_{3},y_{4})=(-4,-2,1,4)$.

Given a partition $\lambda\in \bY$, we consider the following function of $z$:
\begin{equation*}
    G_{\lambda}(z)=\frac{\prod_{i=1}^{n-1}(z-y_{i})}{\prod_{i=1}^{n}(z-x_{i})},
\end{equation*}
where $(x_{1},\dots, x_{n})$ and $(y_{1},\dots, y_{n-1})$ are the local maxima and local minima of the profile of $\lambda$, respectively.
Due to the interlacing property (\ref{eq:interlacing_xy}),
there exists a unique probability measure $\nu_{\lambda}$ on $\bR$ supported on $\{x_{1},\dots, x_{n}\}$ so that the function $G_{\lambda}$ is its Cauchy transform:
\begin{equation*}
    G_{\lambda}(z)=\int_{\bR}\frac{1}{z-x}\nu_{\lambda}(dx).
\end{equation*}
This probability measure $\nu_{\lambda}$ is called the transition measure of $\lambda$.
Obviously, when we expand the Cauchy transform around infinity, the coefficients give the moments of the probability measure in such a way that
\begin{equation*}
    G_{\lambda}(z)=z^{-1}+\sum_{k=1}^{\infty}M_{k}(\lambda)z^{-k-1},\quad \mbox{where}\quad M_{k}(\lambda)=\int_{\bR}x^{k}\nu_{\lambda}(dx),\quad k\in\bZ_{>0}.
\end{equation*}
Note that the property (\ref{eq:sums_xy}) ensures that $M_{1}(\lambda)=0$, i.e., the transition measure $\nu_{\lambda}$ is centered.

We introduce another function; the multiplicative inverse of $G_{\lambda}(z)$:
\begin{equation*}
    H_{\lambda}(z)=\frac{1}{G_{\lambda}(z)}=\frac{\prod_{i=1}^{n}(z-x_{i})}{\prod_{i=1}^{n-1}(z-y_{i})}.
\end{equation*}
The Boolean cumulants $B_{k}(\lambda)$, $k\in \bZ_{>0}$ of $\lambda$ are defined~\cite{SpeicherWoroudi1993} as the coefficients of the expansion of $H_{\lambda}(z)$ around $z=\infty$ so that
\begin{equation*}
    H_{\lambda}(z)=z-\sum_{k=1}^{\infty}B_{k}(\lambda)z^{-k+1}.
\end{equation*}
It is immediate that $B_{1}(\lambda)=M_{1}(\lambda)=0$.

For our purpose, the twisted Boolean cumulants are also important.
They are just the sign change of the Boolean cumulants~\cite{RS2008}: $\hat{B}_{k}(\lambda)=-B_{k}(\lambda)$, $k\in \bZ_{>0}$,
or in other words, the function $H_{\lambda}$ is expanded so that
\begin{equation*}
    H_{\lambda}(z)=z+\sum_{k=1}^{\infty}\hat{B}_{k}(\lambda)z^{-k+1}
\end{equation*}
around $z=\infty$.

\subsection{Main results}
In the previous two subsections, we introduced several families of functions on $\bY$.
One is that of normalized characters $\Sigma_{\pi}$ of the symmetric groups labelled by partitions $\pi\in \cP$.
We also saw the moments, Boolean cumulants, and twisted Boolean cumulants depending on Young diagrams, hence we can regard them as functions on $\bY$.
In the following, borrowing terminology from integrable probability, we often call functions on $\bY$ {\it observables}.

Let us, in particular, take the twisted Boolean cumulants $\{\hat{B}_{k}:k\in \bZ_{>0}\}$.
Our first main result states that any normalized character up to sign is a polynomial of twisted Boolean cumulants with non-negative integer coefficients. Note that the first twisted Boolean cumulant $\hat{B}_{1}\equiv 0$ is a constant function and can be eliminated in the following discussion.

To state the main result, the following notations are convenient: for $\pi\in \cP$, we write $\ell (\pi)$ for the length of $\pi$, i.e., the number of parts in $\pi$. We also set $|\pi|=n$ if $\pi\in \cP_{n}$.
The following result was presented in \cite{RS2008} as a conjecture.
\begin{thm}
\label{thm:main1}
Let $\pi\in \cP$.
There exists a polynomial $P_{\pi}(X_{2},X_{3},\dots)$ with non-negative integer coefficients such that
\begin{equation*}
    (-1)^{\ell (\pi)}\Sigma_{\pi}=P_{\pi}(\hat{B}_{2},\hat{B}_{3},\cdots).
\end{equation*}
\end{thm}

We can say more about the polynomials $P_{\pi}$, $\pi\in \cP$.
\begin{thm}
\label{thm:main2}
Let us define degrees of variables as $\deg X_{i}=i-2$, $i=2,3,\dots$.
Then for any $\pi\in \cP$, all nonzero monomials appearing in $P_{\pi}$ are of degrees at most $|\pi|-\ell (\pi)$.
In particular, each polynomial $P_{\pi}$, $\pi\in \cP$ depends on $X_{2},\dots, X_{|\pi|-\ell (\pi)+2}$.
Furthermore, if we define an automorphism on the ring of polynomials by $\iota\colon X_{k}\mapsto (-1)^{k}X_{k}$, $k=2,3,\dots$,
the polynomials $P_{\pi}$, $\pi\in \cP$ behave so that $\iota (P_{\pi})=(-1)^{|\pi|-\ell (\pi)}P_{\pi}$.
\end{thm}

Conversely, when we expand Boolean cumulants in terms of normalized characters, we again observe non-negative integer coefficients as we state below.
For a partition $\pi\in\cP$ without trivial part, i.e., $\pi_{i}\geq 2$ for all $i=1,\dots, \ell (\pi)$, we set
\begin{equation*}
    \bfB_{\pi}=B_{\pi_{1}}\cdots B_{\pi_{\ell (\pi)}}
\end{equation*}
and call it the Boolean cumulant associated to $\pi$.
\begin{thm}
\label{thm:main3}
For any $\pi\in\cP$ without trivial part, the Boolean cumulant $\bfB_{\pi}$ is expanded in terms of normalized characters so that
\begin{equation*}
    \bfB_{\pi}=\sum_{\substack{\pi'\in \cP\\ |\pi'|\leq |\pi|-\ell(\pi)}}m^{\pi}_{\pi'}\, \Sigma_{\pi'},
\end{equation*}
where $m^{\pi}_{\pi'}\in \bZ_{\geq 0}$ for all appearing $\pi'\in \cP$.
Furthermore, $m^{\pi}_{\pi'}=0$ unless $|\pi'|-\ell(\pi')\equiv |\pi|$ mod $2$.
\end{thm}

% \begin{itemize}
% %    \item \cyan{Tell something about the parity: when defining an automorphism $x_{k}\mapsto (-1)^{k}x_{k}$, we have $P_{\pi}\mapsto (-1)^{|\pi|-\ell (\pi)}P_{\pi}$.}
%     \item \cyan{Expanding a Boolean cumulants in terms of normalized characters, the coefficients are non-negative integers.}
% \end{itemize}

\subsection*{Organization}
In the next Section~\ref{sect:Khovanov_Heisenberg_category}, we introduce Khovanov's Heisenberg category and see some of its properties that are used in subsequent sections.
In Section~\ref{sect:center}, we study the center of Khovanov's Heisenberg category.
In particular, we review the construction of observables on the Young diagrams out of the center of Khovanov's Heisenberg category given by~\cite{KLM2019}.
%In Sect.~\ref{sect:trace}, we briefly recall the trace of Khovanov's Heisenberg category and its action on the center, which will be needed in a technical part of the proof of our result.
We prove Theorems~\ref{thm:main1} and~\ref{thm:main2} in Section~\ref{sect:proofs_main1_main2},
and Theorem~\ref{thm:main3} in Section~\ref{sect:proof_main3}.
Finally, in Section~\ref{sect:discussion}, we make a few comments on related topics and future directions.

\subsection*{Acknowledgments}
This work was supported by Grant-in-Aid for Japan Society for the Promotion of Science Fellows (No.~19J01279) and postdoctoral researcher funding from Academy of Finland (No.~248 130).
The author thanks the anonymous referees for helping the author improve the manuscript, especially by pointing out errors in the first submission.

\section{Khovanov's Heisenberg category}
\label{sect:Khovanov_Heisenberg_category}
\subsection{Definition}
We recall the category $\cH$ that was introduced in \cite{Khovanov2014}.
It is a $\bC$-linear monoidal category generated by objects $Q_{+}$ and $Q_{-}$.
For a sequence $\varepsilon=(\varepsilon_{1},\dots,\varepsilon_{n})$ with $\varepsilon_{i}\in \{+,-\}$, $i=1,\dots, n$, we write $Q_{\varepsilon}=Q_{\varepsilon_{1}}\otimes\cdots\otimes Q_{\varepsilon_{n}}$.
Then objects of $\cH$ are finite direct sums of $Q_{\varepsilon}$ with various sequences $\varepsilon$.
The unit object $\bm{1}=Q_{\emptyset}$ corresponds to the empty sequence.

We will express the objects $Q_{+}$ and $Q_{-}$ by upward and downward arrows as
\begin{equation*}
    Q_{+}=\quad \hackcenter{
    \begin{tikzpicture}
    \draw[thick,->] (0,0)--(0,1);
    \end{tikzpicture}
    }\quad ,\qquad
    Q_{-}=\quad\hackcenter{
    \begin{tikzpicture}
    \draw[thick,->] (0,1)--(0,0);
    \end{tikzpicture}
    }\quad .
\end{equation*}
Then for a sequence $\varepsilon=(\varepsilon_{1},\dots,\varepsilon_{n})$,
the object $Q_{\varepsilon}$ can be depicted by
the corresponding sequence of upward and downward arrows.
For instance, the sequence $\varepsilon=(++-+--)$ gives
\begin{equation*}
    Q_{++-+--}=\quad
    \hackcenter{
    \begin{tikzpicture}
    \draw[thick,->] (0,0)--(0,1);
    \draw[thick,->] (1,0)--(1,1);
    \draw[thick,->] (2,1)--(2,0);
    \draw[thick,->] (3,0)--(3,1);
    \draw[thick,->] (4,1)--(4,0);
    \draw[thick,->] (5,1)--(5,0);
    \end{tikzpicture}
    }\quad .
\end{equation*}

For $\varepsilon=(\varepsilon_{1},\dots,\varepsilon_{m})$ and $\pr{\varepsilon}=(\pr{\varepsilon}_{1},\dots,\pr{\varepsilon}_{n})$,
the morphism space $\Hom_{\cH}(Q_{\varepsilon},Q_{\pr{\varepsilon}})$ is described as follows.
It is a $\bC$-vector space spanned by collections of oriented one-dimensional manifolds immersed in $\bR\times [0,1]$ with endpoints on $\{1,\dots,m\}\times \{0\}$ and $\{1,\dots,n\}\times \{1\}$.
At the endpoints, the orientation of each one-dimensional manifold is indicated by $\varepsilon$ and $\varepsilon'$: at $\{i\}\times\{0\}$, the manifold connected to it is outgoing if $\varepsilon_{i}=+$ and incoming if $\varepsilon_{i}=-$ for all $i=1,\dots, m$. At $\{j\}\times\{1\}$, the manifold is incoming if $\varepsilon'_{j}=+$ and outgoing if $\varepsilon'_{j}=-$ for $j=1,\dots, n$.
% the $\bC$-vector space spanned by oriented chord diagrams that connect the graphical expressions of $Q_{\varepsilon}$ and $Q_{\pr{\varepsilon}}$ with matched orientations, modulo isotopy.
For instance, the following diagram gives a morphism from $Q_{++-+--}$ to $Q_{-+}$:
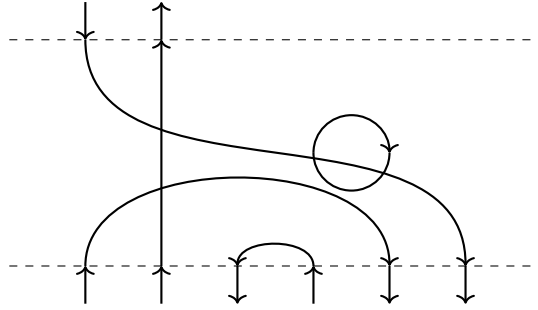
\begin{equation*}
\begin{tikzpicture}
\draw [thick,->] (0,0) to [out=90,in=90] (4,0);
\draw [thick,->] (3,0) to [out=90,in=90] (2,0);
\draw [thick,->] (1,0) to [out=90,in=-90] (1,3);
\draw [thick,->] (0,3) to [out=-90,in=90] (5,0);
\draw [thick,->] (3,1.5) arc (180:0:.5);
\draw [thick] (3,1.5) arc (180:360:.5);
\draw [dashed] (-1,0) -- (6,0);
\draw [dashed] (-1,3) -- (6,3);
\draw [thick,->] (0,-.5) -- (0,0);
\draw [thick,->] (1,-.5) -- (1,0);
\draw [thick,<-] (2,-.5) -- (2,0);
\draw [thick,->] (3,-.5) -- (3,0);
\draw [thick,<-] (4,-.5) -- (4,0);
\draw [thick,<-] (5,-.5) -- (5,0);
\draw [thick,<-] (0,3) -- (0,3.5);
\draw [thick,->] (1,3) -- (1,3.5);
\end{tikzpicture}
\end{equation*}
The morphisms are considered up to isotopy of $\bR\times [0,1]$ preserving the boundaries and are subject to the following local relations:
\begin{equation}
\label{eq:hom_relations_symmetric_group}
\hackcenter{\begin{tikzpicture}
    \draw [thick] (0,0) to [out=90,in=-90] (1,1);
    \draw [thick,->] (1,1) to [out=90,in=-90] (0,2);
    \draw [thick] (1,0) to [out=90,in=-90] (0,1);
    \draw [thick,->] (0,1) to [out=90,in=-90] (1,2);
    \end{tikzpicture}}\quad =\quad
    \hackcenter{
    \begin{tikzpicture}
    \draw [thick,->] (0,0)--(0,2);
%    \draw [thick] (0,1)--(0,2);
    \draw [thick,->] (1,0)--(1,2);
%    \draw [thick] (1,1)--(1,2);
    \end{tikzpicture}
    }\quad ,\qquad\quad
    \hackcenter{
    \begin{tikzpicture}
    \draw [thick,->] (0,0) to [out=90,in=-90] (1,2);
    \draw [thick,->] (1,0) to [out=90,in=-90] (0,2);
    \draw [thick] (.5,0) to [out=90,in=-90] (0,1);
    \draw [thick,->] (0,1) to [out=90,in=-90] (.5,2);
    \end{tikzpicture}
    }\quad =\quad
    \hackcenter{
    \begin{tikzpicture}
    \draw [thick,->] (0,0) to [out=90,in=-90] (1,2);
    \draw [thick,->] (1,0) to [out=90,in=-90] (0,2);
    \draw [thick] (.5,0) to [out=90,in=-90] (1,1);
    \draw [thick,->] (1,1) to [out=90,in=-90] (.5,2);
    \end{tikzpicture}
    }\quad ,
\end{equation}
\vspace{10pt}
\begin{equation}
\label{eq:hom_relations_res_ind}
    \hackcenter{\begin{tikzpicture}
    \draw [thick] (0,0) to [out=90,in=-90] (1,1);
    \draw [thick,->] (1,1) to [out=90,in=-90] (0,2);
    \draw [thick,<-] (1,0) to [out=90,in=-90] (0,1);
    \draw [thick] (0,1) to [out=90,in=-90] (1,2);
    \end{tikzpicture}}\quad =\quad
    \hackcenter{
    \begin{tikzpicture}
    \draw [thick,->] (0,0)--(0,2);
    \draw [thick,<-] (1,0)--(1,2);
    \end{tikzpicture}
    }\quad , \qquad\quad
    \hackcenter{
    \begin{tikzpicture}
    \draw [thick,<-] (0,0) to [out=90,in=-90] (1,1);
    \draw [thick] (1,1) to [out=90,in=-90] (0,2);
    \draw [thick] (1,0) to [out=90,in=-90] (0,1);
    \draw [thick,->] (0,1) to [out=90,in=-90] (1,2);
    \end{tikzpicture}
    }\quad =\quad
    \hackcenter{
    \begin{tikzpicture}
    \draw [thick,<-] (0,0)--(0,2);
    \draw [thick,->] (1,0)--(1,2);
    \end{tikzpicture}
    }\quad -\quad
    \hackcenter{
    \begin{tikzpicture}
    \draw [thick,<-] (0,0) arc (180:0:.5);
    \draw [thick,->] (0,2) arc (180:360:.5);
    \end{tikzpicture}
    }\quad ,
\end{equation}
% \vspace{10pt}
% \begin{equation*}
%     \hackcenter{
%     \begin{tikzpicture}
%     \draw [thick,<-] (0,0) to [out=90,in=-90] (1,1);
%     \draw [thick] (1,1) to [out=90,in=-90] (0,2);
%     \draw [thick] (1,0) to [out=90,in=-90] (0,1);
%     \draw [thick,->] (0,1) to [out=90,in=-90] (1,2);
%     \end{tikzpicture}
%     }\quad =\quad
%     \hackcenter{
%     \begin{tikzpicture}
%     \draw [thick,<-] (0,0)--(0,2);
% %    \draw [thick,<-] (0,1)--(0,2);
%     \draw [thick,->] (1,0)--(1,2);
% %    \draw [thick] (1,1)--(1,2);
%     \end{tikzpicture}
%     }\quad -\quad
%     \hackcenter{
%     \begin{tikzpicture}
%     \draw [thick,<-] (0,0) arc (180:0:.5);
% %    \draw [thick,->] (1,0) arc (0:90:.5);
%     \draw [thick,->] (0,2) arc (180:360:.5);
% %    \draw [thick] (1,2) arc (360:270:.5);
%     \end{tikzpicture}
%     }\quad ,
% \end{equation*}
% \vspace{10pt}
% \begin{equation*}
%     \hackcenter{
%     \begin{tikzpicture}
%     \draw [thick,->] (0,0) to [out=90,in=-90] (1,2);
%     \draw [thick,->] (1,0) to [out=90,in=-90] (0,2);
%     \draw [thick] (.5,0) to [out=90,in=-90] (0,1);
%     \draw [thick,->] (0,1) to [out=90,in=-90] (.5,2);
%     \end{tikzpicture}
%     }\quad =\quad
%     \hackcenter{
%     \begin{tikzpicture}
%     \draw [thick,->] (0,0) to [out=90,in=-90] (1,2);
%     \draw [thick,->] (1,0) to [out=90,in=-90] (0,2);
%     \draw [thick] (.5,0) to [out=90,in=-90] (1,1);
%     \draw [thick,->] (1,1) to [out=90,in=-90] (.5,2);
%     \end{tikzpicture}
%     }\quad ,
% \end{equation*}
\vspace{10pt}
\begin{equation}
\label{eq:hom_relations_circle_turn}
    \hackcenter{
    \begin{tikzpicture}
    \draw [thick,->] (0,0) arc (0:180:.5);
    \draw [thick] (-1,0) arc (180:360:.5);
    \end{tikzpicture}
    }\quad = 1,\qquad\qquad
    \hackcenter{
    \begin{tikzpicture}
    \draw [thick] (0,0) .. controls (0,.5) and (.7,.5) .. (.9,0);
    \draw [thick] (0,0) .. controls (0,-.5) and (.7,-.5) .. (.9,0);
    \draw [thick] (1,-1) .. controls (1,-.4) .. (.9,0);
    \draw [thick,->] (.9,0) .. controls (1,.4) .. (1,1);
    \end{tikzpicture}
    }\quad =0.
\end{equation}

The composition of morphisms is given by the concatenation of diagrams followed by scaling along the vertical direction so that the result fits into $\bR\times [0,1]$.

\begin{rem}
In \cite{Khovanov2014}, studied was the Karoubi envelop of $\cH$, which is usually called Khovanov's Heisenberg category.
In fact, it gives a categorification of the infinite-dimensional Heisenberg algebra.
In the present paper, we do not need the Karoubi envelope since, in the end, we will focus on the center of the category that does not differentiate between $\cH$ and its Karoubi envelope.
Because of this reason, we call $\cH$ Khovanov's Heisenberg category by abuse of terminology.
\end{rem}

\subsection{Some properties}
Here we collect some properties of $\cH$ that we will use in the subsequent parts.

First, we state a trivial fact easily deduced from the definition of $\cH$.
\begin{prop}
\label{prop:symmetric_group_algebra_Hom_space}
For each $n\in \bZ_{>0}$, we write $Q_{(n)}$ for the object corresponding to the sequence $(n)=(\underbrace{+\cdots +}_{n})$.
There exists a unique algebra homomorphism
\begin{equation*}
    \varphi_{n}\colon \bC [S_{n}]\to \End_{\cH}\left(Q_{(n)}\right)
\end{equation*}
such that
\begin{equation*}
    s_{i}\mapsto \quad
    \hackcenter{
    \begin{tikzpicture}
    \draw [thick,->] (0,-.5)--(0,.5);
    \draw (.5,0) node{$\cdots$};
    \draw [thick,->] (1,-.5)--(1,.5);
    \draw [thick,->] (1.5,-.5) to [out=90,in=-90] (2,.5);
    \draw [thick,->] (2,-.5) to [out=90,in=-90] (1.5,.5);
    \draw [thick,->] (2.5,-.5)--(2.5,.5);
    \draw (3,0) node{$\cdots$};
    \draw [thick,->] (3.5,-.5)--(3.5,.5);
    \draw (1.5,-.5) node[below]{$_{i}$};
    \draw (2,-.5) node[below]{$_{i+1}$};
    \draw (0,.7) node[above]{\,};
    \end{tikzpicture}
    }\, ,\quad i=1,\dots, n-1,
\end{equation*}
where $s_{i}=(i\, i+1)$, $i=1,\dots, n-1$ are the standard transpositions.
Similarly, the object $Q_{(-n)}$ corresponds to the sequence $(-n)=(\underbrace{-\cdots -}_{n})$, and there exists a unique algebra homomorphism
\begin{equation*}
    \varphi_{-n}\colon \bC [S_{n}]\to \End_{\cH}\left(Q_{(-n)}\right)
\end{equation*}
characterized by
\begin{equation*}
    s_{i}\mapsto \quad
    \hackcenter{
    \begin{tikzpicture}
    \draw [thick,<-] (0,-.5)--(0,.5);
    \draw (.5,0) node{$\cdots$};
    \draw [thick,<-] (1,-.5)--(1,.5);
    \draw [thick,<-] (1.5,-.5) to [out=90,in=-90] (2,.5);
    \draw [thick,<-] (2,-.5) to [out=90,in=-90] (1.5,.5);
    \draw [thick,<-] (2.5,-.5)--(2.5,.5);
    \draw (3,0) node{$\cdots$};
    \draw [thick,<-] (3.5,-.5)--(3.5,.5);
    \draw (1.5,-.5) node[below]{$_{i}$};
    \draw (2,-.5) node[below]{$_{i+1}$};
    \draw (0,.7) node[above]{\,};
    \end{tikzpicture}
    }\, ,\quad i=1,\dots, n-1.
\end{equation*}
\end{prop}
\begin{proof}
Obvious from the relations (\ref{eq:hom_relations_symmetric_group}).
\end{proof}

For each $\sigma\in S_{n}$, we draw its image under $\varphi_{\pm n}$ as
\begin{equation*}
    \varphi_{n}(\sigma)=
    \hackcenter{
    \begin{tikzpicture}
    \draw[thick,->] (-.75,-.5)--(-.75,.5);
    \draw[thick,->] (.75,-.5)--(.75,.5);
    \fill[gray] (-1,-.25) rectangle (1,.25);
    \draw[thick] (-1,-.25) rectangle (1,.25);
    \draw (0,0) node{$\sigma$};
    \draw (0,-.4) node{$\cdots$};
    \end{tikzpicture}
    }\, ,\quad
    \varphi_{-}(\sigma)=
    \hackcenter{
    \begin{tikzpicture}
    \draw[thick,<-] (-.75,-.5)--(-.75,.5);
    \draw[thick,<-] (.75,-.5)--(.75,.5);
    \fill[gray] (-1,-.25) rectangle (1,.25);
    \draw[thick] (-1,-.25) rectangle (1,.25);
    \draw (0,0) node{$\sigma$};
    \draw (0,-.4) node{$\cdots$};
    \end{tikzpicture}
    }\, .
\end{equation*}

Next, we introduce a shorthand symbol
\begin{equation*}
    \hackcenter{
    \begin{tikzpicture}
    \draw [thick,->] (0,-1)--(0,1);
    \fill[black] (0,0) circle (.1);
    \end{tikzpicture}
    }\quad :=\quad
    \hackcenter{
    \begin{tikzpicture}
    \draw [thick] (.1,0) .. controls (.3,.5) and (1,.5) .. (1,0);
    \draw [thick] (.1,0) .. controls (.3,-.5) and (1,-.5) .. (1,0);
    \draw [thick] (0,-1) .. controls (0,-.4) ..(.1,0);
    \draw [thick,->] (.1,0) .. controls (0,.4) .. (0,1);
    \end{tikzpicture}
    }
\end{equation*}
for the right turn.

\begin{lem}
\label{lem:dot_moving_over_crossing}
We have the following local relations:
\begin{equation*}
\hackcenter{
\begin{tikzpicture}
\draw [thick] (0,0)--(0,.3);
\draw [thick] (1,0)--(1,.3);
\draw [thick] (0,.3) to [out=90,in=-90] (1,1.7);
\draw [thick] (1,.3) to [out=90,in=-90] (0,1.7);
\draw [thick,->] (0,1.7)--(0,2);
\draw [thick,->] (1,1.7)--(1,2);
\fill[black] (0,.3) circle (.1);
\end{tikzpicture}
}\quad=\quad
\hackcenter{
\begin{tikzpicture}
\draw [thick] (0,0)--(0,.3);
\draw [thick] (1,0)--(1,.3);
\draw [thick] (0,.3) to [out=90,in=-90] (1,1.7);
\draw [thick] (1,.3) to [out=90,in=-90] (0,1.7);
\draw [thick,->] (0,1.7)--(0,2);
\draw [thick,->] (1,1.7)--(1,2);
\fill[black] (1,1.7) circle (.1);
\end{tikzpicture}
}\quad +\quad
\hackcenter{
\begin{tikzpicture}
\draw [thick,->] (0,0)--(0,2);
\draw [thick,->] (1,0)--(1,2);
\end{tikzpicture}
}\quad,\qquad\quad
\hackcenter{
\begin{tikzpicture}
\draw [thick] (0,0)--(0,.3);
\draw [thick] (1,0)--(1,.3);
\draw [thick] (0,.3) to [out=90,in=-90] (1,1.7);
\draw [thick] (1,.3) to [out=90,in=-90] (0,1.7);
\draw [thick,->] (0,1.7)--(0,2);
\draw [thick,->] (1,1.7)--(1,2);
\fill[black] (0,1.7) circle (.1);
\end{tikzpicture}
}\quad=\quad
\hackcenter{
\begin{tikzpicture}
\draw [thick] (0,0)--(0,.3);
\draw [thick] (1,0)--(1,.3);
\draw [thick] (0,.3) to [out=90,in=-90] (1,1.7);
\draw [thick] (1,.3) to [out=90,in=-90] (0,1.7);
\draw [thick,->] (0,1.7)--(0,2);
\draw [thick,->] (1,1.7)--(1,2);
\fill[black] (1,.3) circle (.1);
\end{tikzpicture}
}\quad +\quad
\hackcenter{
\begin{tikzpicture}
\draw [thick,->] (0,0)--(0,2);
\draw [thick,->] (1,0)--(1,2);
\end{tikzpicture}
}\quad .
\end{equation*}
\end{lem}
\begin{proof}
As a warm-up for manipulation of diagrams, we show a proof of the first relation.
The second relation in (\ref{eq:hom_relations_res_ind}) gives
\begin{equation*}
    \hackcenter{
\begin{tikzpicture}
\draw [thick] (0,0)--(0,.3);
\draw [thick] (1,0)--(1,.3);
\draw [thick] (0,.3) to [out=90,in=-90] (1,1.7);
\draw [thick] (1,.3) to [out=90,in=-90] (0,1.7);
\draw [thick,->] (0,1.7)--(0,2);
\draw [thick,->] (1,1.7)--(1,2);
\fill[black] (0,.3) circle (.1);
\end{tikzpicture}
}\quad=\quad
    \hackcenter{\begin{tikzpicture}
    \draw [thick,->] (.05,0) ..controls (.2,.5) and (.7,.5) .. (.7,0);
    \draw [thick] (.05,0) ..controls (.2,-.5) and (.7,-.5) .. (.7,0);
    \draw [thick] (0,-.5) .. controls (0,-.3) .. (.05,0);
    \draw [thick] (0,.5) .. controls (0,.3) .. (.05,0);
    \draw [thick,->] (1,-.5) -- (1,0);
    \draw [thick] (1,0)--(1,.5);
    \draw [thick,->] (0,.5) to [out=90,in=-90] (1,2);
    \draw [thick,->] (1,.5) to [out=90,in=-90] (0,2);
    \end{tikzpicture}}\quad=\quad
    \hackcenter{
    \begin{tikzpicture}
    \draw [thick,->] (.05,0) ..controls (.2,.5) and (1,.5) .. (1,0);
    \draw [thick] (.05,0) ..controls (.2,-.5) and (1,-.5) .. (1,0);
    \draw [thick] (0,-.5) .. controls (0,-.3) .. (.05,0);
    \draw [thick] (0,.5) .. controls (0,.3) .. (.05,0);
    \draw [thick,->] (1,-.5) ..controls (.7,-.3) .. (.7,0);
    \draw [thick] (.7,0)--(.7,.5);
    \draw [thick,->] (0,.5) to [out=90,in=-90] (1,2);
    \draw [thick,->] (.7,.5) to [out=90,in=-90] (0,2);
    \end{tikzpicture}
    }\quad +\quad 
    \hackcenter{\begin{tikzpicture}
    \draw [thick] (0,-.5) to [out=90,in=-90] (1,.5);
    \draw [thick] (1,-.5) to [out=90,in=-90] (0,.5);
    \draw [thick,->] (0,.5) to [out=90,in=-90] (1,2);
    \draw [thick,->] (1,.5) to [out=90,in=-90] (0,2);
    \end{tikzpicture}}\quad .
\end{equation*}
We can use (\ref{eq:hom_relations_symmetric_group}) to both terms in the right hand side.
For the first term, we can see that 
\begin{equation*}
    \hackcenter{
    \begin{tikzpicture}
    \draw [thick,->] (.05,0) ..controls (.2,.5) and (1,.5) .. (1,0);
    \draw [thick] (.05,0) ..controls (.2,-.5) and (1,-.5) .. (1,0);
    \draw [thick] (0,-.5) .. controls (0,-.3) .. (.05,0);
    \draw [thick] (0,.5) .. controls (0,.3) .. (.05,0);
    \draw [thick,->] (1,-.5) ..controls (.7,-.3) .. (.7,0);
    \draw [thick] (.7,0)--(.7,.5);
    \draw [thick,->] (0,.5) to [out=90,in=-90] (1,2);
    \draw [thick,->] (.7,.5) to [out=90,in=-90] (0,2);
    \end{tikzpicture}
    }\quad=\quad
\hackcenter{
\begin{tikzpicture}
\draw [thick,->] (.7,-.5) ..controls (1.2, -.1) and (1.5,-.1) .. (1.5,-.5);
\draw [thick] (.7,-.5) ..controls (1.2, -.9) and (1.5,-.9) .. (1.5,-.5);
\draw [thick] (0,-1.25) ..controls (0,-1) .. (.7,-.5);
\draw [thick,->] (1,-1.25) to [out=90,in=-90] (0,-.5);
\draw [thick] (0,-.5) to [out=90,in=-90] (.7,.5);
\draw [thick] (.7,-.5) ..controls (0,.1) .. (0,.5);
\draw [thick,->] (.7,.5) to [out=90,in=-90] (0,1.25);
\draw [thick,->] (0,.5) to [out=90,in=-90] (1,1.25);
\end{tikzpicture}
}\quad =\quad
\hackcenter{
\begin{tikzpicture}
\draw [thick] (0,0)--(0,.3);
\draw [thick] (1,0)--(1,.3);
\draw [thick] (0,.3) to [out=90,in=-90] (1,1.7);
\draw [thick] (1,.3) to [out=90,in=-90] (0,1.7);
\draw [thick,->] (0,1.7)--(0,2);
\draw [thick,->] (1,1.7)--(1,2);
\fill[black] (1,1.7) circle (.1);
\end{tikzpicture}
}\quad ,
\end{equation*}
where we used the second relation in (\ref{eq:hom_relations_symmetric_group}) in the former equality to ``pull the loop rightward",
and the first relation in the latter equality to resolve the consecutive two crossings.
The latter term is obviously
\begin{equation*}
    \hackcenter{\begin{tikzpicture}
    \draw [thick] (0,-.5) to [out=90,in=-90] (1,.5);
    \draw [thick] (1,-.5) to [out=90,in=-90] (0,.5);
    \draw [thick,->] (0,.5) to [out=90,in=-90] (1,2);
    \draw [thick,->] (1,.5) to [out=90,in=-90] (0,2);
    \end{tikzpicture}}\quad =\quad
    \hackcenter{
    \begin{tikzpicture}
    \draw [thick,->] (0,-.5) -- (0,2);
    \draw [thick,->] (1,-.5) -- (1,2);
    \end{tikzpicture}
    }\quad .
\end{equation*}
Hence, we obtain the desired identity.
\end{proof}

We also introduce a notation for multiple insertion of dots:
\begin{equation*}
    \hackcenter{
    \begin{tikzpicture}
    \draw [thick,->] (0,-1)--(0,1);
    \fill[black] (0,0) circle (.1);
    \draw (0,0) node[left]{$k$};
    \end{tikzpicture}
    }\quad :=\quad
    \hackcenter{
    \begin{tikzpicture}
    \draw [thick] (0,-1)--(0,-.3);
    \draw [dotted] (0,-.3)--(0,.3);
    \draw [thick,->] (0,.3)--(0,1);
    \fill[black] (0,-.6) circle (.1);
    \fill[black] (0,.6) circle (.1);
    \draw [semithick](.2,.6) to [out=0,in=180] (.42,0);
    \draw [semithick](.2,-.6) to [out=0,in=180] (.42,0);
    \draw (.42,0) node[right]{$k$};
    \end{tikzpicture}
    }\quad ,\qquad k\in \bZ_{\geq 0}.
\end{equation*}

\begin{lem}
\label{lem:moving_dotted_circle_one_step}
For any $k\in \bZ_{\geq 0}$,
we have the following relation:
\begin{align*}
    \hackcenter{
    \begin{tikzpicture}
    \draw [thick,->] (0,-1)--(0,1);
    \draw [thick,->] (.7,0) arc (180:0:.5);
    \draw [thick] (1.7,0) arc (360:180:.5);
    \fill [black] (.7,0) circle (.1);
    \draw (.7,0) node[left]{$k$};
    \end{tikzpicture}
    }\quad =\quad
    \hackcenter{
    \begin{tikzpicture}
    \draw [thick,->] (2,-1)--(2,1);
    \draw [thick,->] (.5,0) arc (180:0:.5);
    \draw [thick] (1.5,0) arc (360:180:.5);
    \fill [black] (.5,0) circle (.1);
    \draw (.5,0) node[left]{$k$};
    \end{tikzpicture}
    }\quad -(k+1)\,
    \hackcenter{
    \begin{tikzpicture}
    \draw [thick,->] (0,-1)--(0,1);
    \fill[black] (0,0) circle (.1);
    \draw (0,0) node[left]{$k$};
    \end{tikzpicture}
    }\quad +\sum_{b=0}^{k-2}(b+1)\,
    \hackcenter{
    \begin{tikzpicture}
    \draw [thick,->] (0,-1)--(0,1);
    \draw [thick,->] (.7,0) arc (180:0:.5);
    \draw [thick] (1.7,0) arc (360:180:.5);
    \fill [black] (1.2,-.5) circle (.1);
    \draw (1.2,-.5) node[below]{$k-2-b$};
    \fill [black] (0,0) circle (.1);
    \draw (0,0) node[left]{$b$};
    \end{tikzpicture}
    }\, ,
\end{align*}
where the last sum is understood as empty unless $k\geq 2$.
\end{lem}
\begin{proof}
For each $i=0,1,\dots, k-1$, let us consider the following diagram:
\begin{equation*}
    \hackcenter{\begin{tikzpicture}
    \draw [thick,->] (0,-1)--(0,1);
    \draw [thick,->] (-.3,0) arc (180:0:.5);
    \draw [thick] (.7,0) arc (360:180:.5);
    \fill [black] (-.3,0) circle (.1);
    \fill [black] (.2,-.5) circle (.1);
    \draw (-.3,0) node[left]{$i$};
    \draw (.2,-.5) node[below right]{$k-i$};
    \end{tikzpicture}}\, .
\end{equation*}
Moving one dot from the right of the upward arrow to the left using Lemma~\ref{lem:dot_moving_over_crossing}, it becomes
\begin{equation*}
    \hackcenter{\begin{tikzpicture}
    \draw [thick,->] (0,-1)--(0,1);
    \draw [thick,->] (-.3,0) arc (180:0:.5);
    \draw [thick] (.7,0) arc (360:180:.5);
    \fill [black] (-.3,0) circle (.1);
    \fill [black] (.2,-.5) circle (.1);
    \draw (-.3,0) node[left]{$i$};
    \draw (.2,-.5) node[below right]{$k-i$};
    \end{tikzpicture}}\quad =\quad
    \hackcenter{\begin{tikzpicture}
    \draw [thick,->] (0,-1)--(0,1);
    \draw [thick,->] (-.3,0) arc (180:0:.5);
    \draw [thick] (.7,0) arc (360:180:.5);
    \fill [black] (-.3,0) circle (.1);
    \fill [black] (.2,-.5) circle (.1);
    \draw (-.3,0) node[left]{$i+1$};
    \draw (.2,-.5) node[below right]{$k-i-1$};
    \end{tikzpicture}}\quad -\quad
    \hackcenter{
    \begin{tikzpicture}
    \draw [thick] (0,0) .. controls (.5,.7) and (1,.5) .. (1,0);
    \draw [thick] (0,0) .. controls (.5,-.7) and (1,-.5) .. (1,0);
    \draw [thick] (-.3,-1) .. controls (-.3,-1) and (-.3,-.4) ..(0,0);
    \draw [thick,->] (0,0) .. controls (-.3,.4) and (-.3,1) ..(-.3,1);
    \fill [black] (1,0) circle (.1);
    \fill [black] (-.22,-.5) circle (.1);
    \draw (1,0) node[right]{$k-i-1$};
    \draw (-.22,-.5) node[left]{$i$};
    \end{tikzpicture}
    }\, .
\end{equation*}
The same Lemma~\ref{lem:dot_moving_over_crossing} allows one to transform the second diagram in the right hand side in the following way:
\begin{equation*}
    \hackcenter{
    \begin{tikzpicture}
    \draw [thick] (0,0) .. controls (.5,.7) and (1,.5) .. (1,0);
    \draw [thick] (0,0) .. controls (.5,-.7) and (1,-.5) .. (1,0);
    \draw [thick] (-.3,-1) .. controls (-.3,-1) and (-.3,-.4) ..(0,0);
    \draw [thick,->] (0,0) .. controls (-.3,.4) and (-.3,1) ..(-.3,1);
    \fill [black] (1,0) circle (.1);
    \fill [black] (-.22,-.5) circle (.1);
    \draw (1,0) node[right]{$k-i-1$};
    \draw (-.22,-.5) node[left]{$i$};
    \end{tikzpicture}
    }\quad =\quad 
    \hackcenter{
    \begin{tikzpicture}
    \draw [thick,->] (0,-1) -- (0,1);
    \fill [black] (0,0) circle (.1);
    \draw (0,0) node[left]{$k$};
    \end{tikzpicture}
    }\quad -\sum_{b=i}^{k-2}\quad
    \hackcenter{
    \begin{tikzpicture}
    \draw [thick,->] (0,-1) -- (0,1);
    \draw [thick,->] (.5,0) arc (180:0:.5);
    \draw [thick] (1.5,0) arc (360:180:.5);
    \fill [black] (0,0) circle (.1);
    \fill [black] (1,-.5) circle (.1);
    \draw (0,0) node[left]{$b$};
    \draw (1,-.5) node[below right]{$k-b-2$};
    \end{tikzpicture}
    }\, .
\end{equation*}
Consequently, we obtain the relation
\begin{equation}
\label{eq:moving_dot_on_circle}
    \hackcenter{\begin{tikzpicture}
    \draw [thick,->] (0,-1)--(0,1);
    \draw [thick,->] (-.3,0) arc (180:0:.5);
    \draw [thick] (.7,0) arc (360:180:.5);
    \fill [black] (-.3,0) circle (.1);
    \fill [black] (.2,-.5) circle (.1);
    \draw (-.3,0) node[left]{$i$};
    \draw (.2,-.5) node[below right]{$k-i$};
    \end{tikzpicture}}\, =\quad
    \hackcenter{\begin{tikzpicture}
    \draw [thick,->] (0,-1)--(0,1);
    \draw [thick,->] (-.3,0) arc (180:0:.5);
    \draw [thick] (.7,0) arc (360:180:.5);
    \fill [black] (-.3,0) circle (.1);
    \fill [black] (.2,-.5) circle (.1);
    \draw (-.3,0) node[left]{$i+1$};
    \draw (.2,-.5) node[below right]{$k-i-1$};
    \end{tikzpicture}} \, -\quad
    \hackcenter{
    \begin{tikzpicture}
    \draw [thick,->] (0,-1) -- (0,1);
    \fill [black] (0,0) circle (.1);
    \draw (0,0) node[left]{$k$};
    \end{tikzpicture}
    }\quad +\sum_{b=i}^{k-2}\quad
    \hackcenter{
    \begin{tikzpicture}
    \draw [thick,->] (0,-1) -- (0,1);
    \draw [thick,->] (.5,0) arc (180:0:.5);
    \draw [thick] (1.5,0) arc (360:180:.5);
    \fill [black] (0,0) circle (.1);
    \fill [black] (1,-.5) circle (.1);
    \draw (0,0) node[left]{$b$};
    \draw (1,-.5) node[below right]{$k-b-2$};
    \end{tikzpicture}
    }\, ,
\end{equation}
where the last sum in the right hand side is understood as empty if $i=k-1$.
When we sum up (\ref{eq:moving_dot_on_circle}) from $i=0$ to $i=k-1$,
we can see that
\begin{align*}
    \hackcenter{
    \begin{tikzpicture}
    \draw [thick,->] (0,-1)--(0,1);
    \draw [thick,->] (.7,0) arc (180:0:.5);
    \draw [thick] (1.7,0) arc (360:180:.5);
    \fill [black] (.7,0) circle (.1);
    \draw (.7,0) node[left]{$k$};
    \end{tikzpicture}
    }\quad &=\quad
    \hackcenter{\begin{tikzpicture}
    \draw [thick,->] (0,-1)--(0,1);
    \draw [thick,->] (-.3,0) arc (180:0:.5);
    \draw [thick] (.7,0) arc (360:180:.5);
    \fill [black] (.2,-.5) circle (.1);
    \draw (.2,-.5) node[below right]{$k$};
    \end{tikzpicture}} \\
    &=\quad 
    \hackcenter{\begin{tikzpicture}
    \draw [thick,->] (0,-1)--(0,1);
    \draw [thick,->] (-.7,0) arc (180:0:.5);
    \draw [thick] (.3,0) arc (360:180:.5);
    \fill [black] (-.7,0) circle (.1);
    \draw (-.7,0) node[left]{$k$};
    \end{tikzpicture}}\quad -k\,
    \hackcenter{
    \begin{tikzpicture}
    \draw [thick,->] (0,-1) -- (0,1);
    \fill [black] (0,0) circle (.1);
    \draw (0,0) node[left]{$k$};
    \end{tikzpicture}
    }\quad +\sum_{b=0}^{k-2}(b+1)\,
    \hackcenter{
    \begin{tikzpicture}
    \draw [thick,->] (0,-1) -- (0,1);
    \draw [thick,->] (.5,0) arc (180:0:.5);
    \draw [thick] (1.5,0) arc (360:180:.5);
    \fill [black] (0,0) circle (.1);
    \fill [black] (1,-.5) circle (.1);
    \draw (0,0) node[left]{$b$};
    \draw (1,-.5) node[below right]{$k-b-2$};
    \end{tikzpicture}
    } \\
    &=\quad
    \hackcenter{
    \begin{tikzpicture}
    \draw [thick,->] (2,-1)--(2,1);
    \draw [thick,->] (.5,0) arc (180:0:.5);
    \draw [thick] (1.5,0) arc (360:180:.5);
    \fill [black] (.5,0) circle (.1);
    \draw (.5,0) node[left]{$k$};
    \end{tikzpicture}
    }\quad -(k+1)\,
    \hackcenter{
    \begin{tikzpicture}
    \draw [thick,->] (0,-1)--(0,1);
    \fill[black] (0,0) circle (.1);
    \draw (0,0) node[left]{$k$};
    \end{tikzpicture}
    }\quad +\sum_{b=0}^{k-2}(b+1)\,
    \hackcenter{
    \begin{tikzpicture}
    \draw [thick,->] (0,-1)--(0,1);
    \draw [thick,->] (.7,0) arc (180:0:.5);
    \draw [thick] (1.7,0) arc (360:180:.5);
    \fill [black] (1.2,-.5) circle (.1);
    \draw (1.2,-.5) node[below]{$k-2-b$};
    \fill [black] (0,0) circle (.1);
    \draw (0,0) node[left]{$b$};
    \end{tikzpicture}
    }
\end{align*}
as is desired.
\end{proof}

\begin{prop}
\label{prop:bubble_move_non-negativity}
For any $k\in \bZ_{\geq 0}$, we have the relation
\begin{equation*}
    \hackcenter{
    \begin{tikzpicture}
    \draw [thick,->] (0,-1)--(0,1);
    \draw [thick,->] (.7,0) arc (180:0:.5);
    \draw [thick] (1.7,0) arc (360:180:.5);
    \fill [black] (.7,0) circle (.1);
    \draw (.7,0) node[left]{$k$};
    \end{tikzpicture}
    }\quad =\sum_{\substack{i,j\in \bZ_{\geq 0}\\i+j\leq k}}m_{ij}\,
    \hackcenter{
    \begin{tikzpicture}
    \draw [thick,->] (2,-1)--(2,1);
    \draw [thick,->] (.5,0) arc (180:0:.5);
    \draw [thick] (1.5,0) arc (360:180:.5);
    \fill [black] (.5,0) circle (.1);
    \fill [black] (2,0) circle (.1);
    \draw (.5,0) node[left]{$i$};
    \draw (2,0) node[right]{$j$};
    \end{tikzpicture}
    }\quad -\sum_{l=0}^{k}n_{l}\,
    \hackcenter{
    \begin{tikzpicture}
    \draw [thick,->] (0,-1)--(0,1);
    \fill[black] (0,0) circle (.1);
    \draw (0,0) node[left]{$l$};
    \end{tikzpicture}
    },
\end{equation*}
where $m_{ij},n_{l}\in \bZ_{\geq 0}$ for all $i,j,l$.
\end{prop}
\begin{proof}
When $k=0$, the assertion is obvious from Lemma~\ref{lem:moving_dotted_circle_one_step}.
Assume that the assertion is true up to some $k$,
then from Lemma~~\ref{lem:moving_dotted_circle_one_step}, we have
\begin{equation*}
    \hackcenter{
    \begin{tikzpicture}
    \draw [thick,->] (0,-1)--(0,1);
    \draw [thick,->] (.5,0) arc (180:0:.5);
    \draw [thick] (1.5,0) arc (360:180:.5);
    \fill [black] (1,-.5) circle (.1);
    \draw (1,-.5) node[below]{$k+1$};
    \end{tikzpicture}}
    \quad =\quad
    \hackcenter{
    \begin{tikzpicture}
    \draw [thick,->] (2,-1)--(2,1);
    \draw [thick,->] (.5,0) arc (180:0:.5);
    \draw [thick] (1.5,0) arc (360:180:.5);
    \fill [black] (1,-.5) circle (.1);
    \draw (1,-.5) node[below]{$k+1$};
    \end{tikzpicture}
    }\quad -(k+2)\,
    \hackcenter{
    \begin{tikzpicture}
    \draw [thick,->] (0,-1)--(0,1);
    \fill[black] (0,0) circle (.1);
    \draw (0,0) node[left]{$k+1$};
    \end{tikzpicture}
    }\quad +\sum_{b=0}^{k-1}(b+1)\,
    \hackcenter{
    \begin{tikzpicture}
    \draw [thick,->] (0,-1)--(0,1);
    \draw [thick,->] (.7,0) arc (180:0:.5);
    \draw [thick] (1.7,0) arc (360:180:.5);
    \fill [black] (1.2,-.5) circle (.1);
    \draw (1.2,-.5) node[below]{$k-1-b$};
    \fill [black] (0,0) circle (.1);
    \draw (0,0) node[left]{$b$};
    \end{tikzpicture}
    }\, .
\end{equation*}
On each term in the last sum, we apply the induction hypothesis in a local way: we can apply the induction hypothesis in the blue dashed box below to get
\begin{align*}
    \hackcenter{
    \begin{tikzpicture}
    \draw [thick,->] (0,-1)--(0,2);
    \draw [thick,->] (.7,0) arc (180:0:.5);
    \draw [thick] (1.7,0) arc (360:180:.5);
    \fill [black] (1.2,-.5) circle (.1);
    \draw (1.2,-.5) node[below]{$k-1-b$};
    \fill [black] (0,1) circle (.1);
    \draw (0,1) node[left]{$b$};
    \draw [thick, blue,dashed] (-.5,.7)--(2.2,.7)--(2.2,-1)--(-.5,-1)--cycle;
    \end{tikzpicture}
    }\quad &=
    \sum_{\substack{i,j\in \bZ_{\geq 0}\\i+j\leq k-1-b}}m^{b}_{i,j}\,
    \hackcenter{
    \begin{tikzpicture}
    \draw [thick,->] (2,-1)--(2,2);
    \draw [thick,->] (.5,0) arc (180:0:.5);
    \draw [thick] (1.5,0) arc (360:180:.5);
    \fill [black] (.5,0) circle (.1);
    \fill [black] (2,0) circle (.1);
    \fill [black] (2,1) circle (.1);
    \draw (.5,0) node[left]{$i$};
    \draw (2,0) node[right]{$j$};
    \draw (2,1) node[right]{$b$};
    \end{tikzpicture}
    }\quad -\sum_{l=0}^{k-1-b}n^{b}_{l}\,
    \hackcenter{
    \begin{tikzpicture}
    \draw [thick,->] (0,-1)--(0,2);
    \fill[black] (0,0) circle (.1);
    \fill[black] (0,1) circle (.1);
    \draw (0,0) node[left]{$l$};
    \draw (0,1) node[left]{$b$};
    \end{tikzpicture}
    } \\
    &=\sum_{\substack{i,j\in \bZ_{\geq 0},\, j\geq b\\i+j\leq k-1}}m^{b}_{i,j-b} \,
    \hackcenter{
    \begin{tikzpicture}
    \draw [thick,->] (2,-1)--(2,1);
    \draw [thick,->] (.5,0) arc (180:0:.5);
    \draw [thick] (1.5,0) arc (360:180:.5);
    \fill [black] (.5,0) circle (.1);
    \fill [black] (2,0) circle (.1);
    \draw (.5,0) node[left]{$i$};
    \draw (2,0) node[right]{$j$};
    \end{tikzpicture}
    }\quad -\sum_{l=b}^{k-1}n^{b}_{l-b}\,
    \hackcenter{
    \begin{tikzpicture}
    \draw [thick,->] (0,-1)--(0,1);
    \fill[black] (0,0) circle (.1);
    \draw (0,0) node[left]{$l$};
    \end{tikzpicture}
    },
\end{align*}
where $m_{ij}^{b},n^{b}_{l}\in \bZ_{\geq 0}$.
Therefore, the assertion is also true for $k+1$.
\end{proof}

\section{Center of Khovanov's Heisenberg category}
\label{sect:center}
\subsection{Basic description}
By definition, the center of the category $\cH$ is the commutative algebra
\begin{equation*}
    Z(\cH)= \End_{\cH}(\bm{1}).
\end{equation*}

Let us recall the generators of $Z(\cH)$ studied in \cite{Khovanov2014}.
We set
\begin{equation*}
    c_{k}:=\,
    \hackcenter{
    \begin{tikzpicture}
    \draw [thick,->] (0,0) arc (180:0:.5);
    \draw [thick] (1,0) arc (360:180:.5);
    \fill [black] (0,0) circle (.1);
    \draw (0,0) node[left]{$k$};
    \end{tikzpicture}
    }\, ,\quad
    \tilde{c}_{k}:=\,
    \hackcenter{
    \begin{tikzpicture}
    \draw [thick,->] (0,0) arc (180:360:.5);
    \draw [thick] (1,0) arc (0:180:.5);
    \fill [black] (0,0) circle (.1);
    \draw (0,0) node[left]{$k$};
    \end{tikzpicture}
    }\, ,\quad k\in \bZ_{\geq 0}.
\end{equation*}
From the first relation in (\ref{eq:hom_relations_circle_turn}), we have $\tilde{c}_{0}=1$.
Also applying the second relation in (\ref{eq:hom_relations_circle_turn}), we can see that $\tilde{c}_{1}=0$.
Hence as for the series $\tilde{c}_{k}$, only those for $k\geq 2$
are nontrivial algebra elements.

\begin{prop}[\cite{Khovanov2014}]
\label{prop:center_generators}
We have isomorphisms
\begin{equation*}
    Z(\cH)\simeq \bC [c_{0},c_{1},\dots]\simeq \bC [\tilde{c}_{2},\tilde{c}_{3},\dots]
\end{equation*}
of commutative algebras.
The two series of generators $\{c_{k}\}_{k\geq 0}$ and $\{\tilde{c}_{k}\}_{k\geq 2}$ satisfy the following relations
\begin{equation*}
    \tilde{c}_{k+1}=\sum_{i=0}^{k-1}\tilde{c}_{i}c_{k-1-i},\quad k\geq 1.
\end{equation*}
\end{prop}

\subsection{Representation of \texorpdfstring{$\cH$}{H}}
\label{subsect:rep_of_H}
In this and next sections, we expose a remarkable description of $Z(\cH)$ found in \cite{KLM2019}, but in a slightly different language.
For that purpose, we first consider a certain representation of the category denoted by $\cH$.

Given an algebra $A$, we write $A$-$\catMod$ for the category of finite rank modules over $A$.
For $n\in \bZ_{\geq 0}$, the induction functor
\begin{align*}
    \funInd_{n}^{n+1}\colon \bC [S_{n}]\mbox{-}\catMod \to \bC [S_{n+1}]\mbox{-}\catMod
\end{align*}
assigns to a $\bC[S_{n}]$-module $M$ the $\bC[S_{n+1}]$-module $\bC[S_{n+1}]\otimes_{\bC[S_{n}]}M$ according to the natural embedding $S_{n}\subset S_{n+1}$ fixing the last letter.
On the other hand, the restriction functor
\begin{equation*}
    \funRes^{n+1}_{n}\colon \bC [S_{n+1}]\mbox{-}\catMod \to \bC [S_{n}]\mbox{-}\catMod
\end{equation*}
sends a $\bC[S_{n+1}]$-module $M$ to the same vector space considered a $\bC[S_{n}]$-module.

The representation category of interest is the direct sum
\begin{equation*}
    \cC:=\bigoplus_{n\geq 0}\bC [S_{n}]\mbox{-}\catMod
\end{equation*}
of module categories.
To be specific, its object is a direct sum $\bigoplus_{n\geq 0}M_{n}$ of $M_{n}\in \bC [S_{n}]\mbox{-}\catMod$, $n\geq 0$ such that $M_{n}=0$ except finitely many $n$. A morphism is also a direct sum $\bigoplus_{n\geq 0}f_{n}$, where each $f_{n}$ is a morphism in the category $\bC [S_{n}]\mbox{-}\catMod$.

By a representation of $\cH$ on $\cC$, we mean a monoidal functor
\begin{equation*}
    \cF\colon \cH\to \cE nd(\cC),
\end{equation*}
where $\cE nd(\cC)$ is the category of endofunctors on $\cC$. Note that the monoidal structure of this category is given by the composition of functors.

At the object level, we define $\cF$ by
\begin{align*}
    \cF (Q_{+}):=\bigoplus_{n\geq 0}\funInd_{n}^{n+1},\quad
    \cF (Q_{-}):=\bigoplus_{n\geq 0}\funRes^{n+1}_{n}
\end{align*}
and extend them by linearity and monoidal structure.
For a morphism $f\colon X\to Y$ in $\cH$, $\cF(f)$ is a natural transformation of functors $\cF(f)\colon \cF(X)\Rightarrow \cF(Y)$.
To describe such a natural transformation, it is convenient to fix $M\in \bC[S_{n}]\mbox{-}\catMod$.
Then, we define $\cF (f)_{M}\colon \cF (X)(M)\to \cF (Y)(M)$ for the following pieces:
\begin{align*}
    \cF \left(
    \hackcenter{
    \begin{tikzpicture}
    \draw [thick,->] (-.4,-.4)--(.4,.4);
    \draw [thick,->] (.4,-.4)--(-.4,.4);
    \end{tikzpicture}
    }
    \right)_{M}&\colon \bC [S_{n+2}]\otimes_{\bC[S_{n}]}M\to \bC [S_{n+2}]\otimes_{\bC[S_{n}]}M;\quad g\otimes v\mapsto gs_{n+1}\otimes v, \\
    \cF \left(
    \hackcenter{
    \begin{tikzpicture}
    \draw [thick,->] (0,0) arc (180:0:.4);
    \end{tikzpicture}
    }
    \right)_{M}&\colon \bC [S_{n}]\otimes_{\bC [S_{n-1}]}M\to M;\quad g\otimes v\mapsto gv, \\
    \cF \left(
    \hackcenter{
    \begin{tikzpicture}
    \draw [thick,->] (0,0) arc (180:360:.4);
    \end{tikzpicture}
    }
    \right)_{M}&\colon M\to \funRes^{n+1}_{n}\left(\bC[S_{n+1}]\otimes_{\bC [S_{n}]}M\right);\quad v\mapsto 1\otimes v,\\
    \cF \left(
    \hackcenter{
    \begin{tikzpicture}
    \draw [thick,->] (1.2,0) arc (0:180:.4);
    \end{tikzpicture}
    }
    \right)_{M}&\colon \funRes^{n+1}_{n}\left(\bC[S_{n+1}]\otimes_{\bC [S_{n}]}M\right)\to M;\quad g\otimes v\mapsto \mathrm{pr}_{n}[g]v, \\
    \cF \left(
    \hackcenter{
    \begin{tikzpicture}
    \draw [thick,->] (1.2,0) arc (360:180:.4);
    \end{tikzpicture}
    }
    \right)_{M}&\colon M\to \bC[S_{n}]\otimes_{\bC[S_{n-1}]}M;\quad v\mapsto\sum_{i=1}^{n+1}s_{i}s_{i+1}\cdots s_{n}\otimes s_{n}\cdots s_{i+1}s_{i}v.
\end{align*}
Here, $\mathrm{pr}_{n}\colon \bC [S_{n+1}]\to \bC [S_{n}]$ is the conditional expectation defined by
\begin{equation*}
    g\mapsto
    \begin{cases}
    g & g\in S_{n},\\
    0 & g\not\in S_{n},
    \end{cases}\quad g\in S_{n+1},
\end{equation*}
This definition of $\cF (f)_{M}$ is naturally extended to an arbitrary morphism $f$ in $\cH$ and an object $M$ in $\cC$.

\begin{prop}[\cite{Khovanov2014}]
The functor $\cF\colon \cH\to \cE nd (\cC)$ gives a representation of the monoidal category $\cF$ on $\cC$.
\end{prop}

In fact, the relations imposed on the morphism space of $\cH$
illustrate the {\it commutation relations} of induction and restriction functors.

% Note that the endomorphism ring $\End_{\cS_{k}}(\funId_{\bC [S_{k}]\mbox{-}\catMod})$ consists of bimodule automorphisms on $\bC [S_{k}]$ over $(\bC [S_{k}],\bC [S_{k}])$.
% In other words, the endomorphism ring is simply the center of $\bC [S_{k}]$:
% \begin{equation*}
%     \End_{\cS_{k}}(\funId_{\bC [S_{k}]\mbox{-}\catMod})\simeq Z(\bC [S_{k}]).
% \end{equation*}
% Therefore the functor $\cF_{k}$ induces the ring homomorphism (denoted by the same symbol)
% \begin{equation*}
%     \cF_{k}\colon Z(\cH)\to Z(\bC [S_{k}]).
% \end{equation*}

\subsection{Construction by Kvinge--Licata--Mitchell}
Applying the functor to the center $Z (\cH)=\End_{\cH}(\bm{1})$, we get an algebra homomorphism
\begin{equation*}
    \cF\colon Z(\cH)\to \End_{\cE nd(\cC)}(\funId).
\end{equation*}
We write $\mathrm{Fun}(\bY)$ for the ring of complex-valued functions on $\bY$.
\begin{lem}
We have an isomorphism of algebras $\End_{\cE nd(\cC)}(\funId)\simeq \mathrm{Fun}(\bY)$,
where $\mathrm{Fun}(\bY)$ is equipped with the point-wise product.
\end{lem}
\begin{proof}
Let $\eta \colon \funId\Rightarrow \funId$ be a natural transformation.
For each $\lambda\in \bY$, there is $f_{\eta}(\lambda)\in\bC$ determined by $\eta_{V^{\lambda}}=f_{\eta}(\lambda)\id_{V^{\lambda}}$.
For two natural transformations $\eta$ and $\nu$, it is obvious that $f_{\eta\circ\nu}(\lambda)=f_{\eta}(\lambda)f_{\nu}(\lambda)$ for all $\lambda\in \bY$.
Therefore, the assignment $\eta \mapsto f_{\eta}$ is an algebraic homomorphism.

Conversely, given a function $f\in \mathrm{Fun}(\bY)$, we can construct a natural transformation $\eta^{f}\colon \funId\Rightarrow\funId$ by $\eta^{f}_{V^{\lambda}}=f(\lambda)\id_{V^{\lambda}}$, $\lambda\in\bY$ and by naturally extending this to arbitrary objects.
This assignment $f\mapsto \eta^{f}$ is also an algebra homomorphism.

It is obvious from construction that the assignments $\eta\mapsto f_{\eta}$ and $f\mapsto \eta^{f}$ are the inverse of each other, proving the desired result.
\end{proof}

Consequently, we get an algebra homomorphism
\begin{equation*}
    \cF\colon Z(\cH)\to \mathrm{Fun}(\bY),
\end{equation*}
which is denoted by the same symbol as before.

% For $\lambda\in \bY_{k}$,
% we can consider the character $\chi^{\lambda}$ as a function on $Z(\bC [S_{k}])$ by linearly extending it to $\bC [S_{k}]$ and restricting it to $Z(\bC [S_{k}])$.
% Furthermore, when we recall the isomorpshisms (see e.g.~\cite{Lam2001})
% \begin{equation*}
%     \bC [S_{k}]\simeq \bigoplus_{\nu\in \bY_{k}}\End (V^{\mu}),\quad 
%     Z(\bC [S_{k}]) \simeq \bigoplus_{\nu\in \bY_{k}}\bC\, \id_{V^{\mu}}
% \end{equation*}
% of algebras, we can see that the irreducible character $\chi^{\lambda}$ restricted on $ Z(\bC [S_{k}])$ is an algebra homomorphism.

% We write $\bC [\bY]$ for the ring of complex-valued functions on $\bY$.
% \begin{prop}[\cite{KLM2019}]
% Define a map $\Phi\colon Z(\cH)\to \bC[\bY]$ by
% \begin{equation*}
%     (\Phi (x))(\lambda):=(\chi^{\lambda}\circ \cF_{k})(x),\quad x\in Z(\cH),\quad \lambda\in \bY_{k},\quad k\in \bZ_{\geq 0}.
% \end{equation*}
% Then it is an injective algebraic homomorphism.
% \end{prop}
% \begin{proof}
% Since the functor $\cF_{k}$ induces an algebraic homomorphism $Z(\cH)\to Z(\bC [S_{k}])$, and the character $\chi^{\lambda}$ is also an algebraic homomorphism on $Z(\bC [S_{k}])$, it is obvious that their composition, $\Phi$ is an algebraic homomorphism (the algebraic structure of $\bC[\bY]$ is defined by the pointwise multiplication.)
% In \cite{KLM2019}, the authors showed that algebraically independent generators of $Z(\cH)$ are sent to algebraically independent elements of $\bC [\bY]$, which implies the injectivity.
% \end{proof}

Let us introduce some notations.
First, we set for each $k\in \bZ_{>0}$,
\vspace{-20pt}
\begin{equation*}
    \hackcenter{
    \begin{tikzpicture}
    \draw [thick,->] (0,-.5)--(0,.5);
    \draw [thick,->] (.5,-.5)--(.5,.5);
    \draw [thick,->] (2,-.5)--(2,.5);
    \fill [white] (-.25,-.25) rectangle (2.25,.25);
    \draw [thick] (-.25,-.25) rectangle (2.25,.25);
    \draw (1,0) node{$k$};
    \draw (1.25,-.4) node{$\cdots$};
    \draw (1.25,.4) node{$\cdots$};
    % \draw [semithick] (-.1,-.55) to [out=-90,in=90] (1,-1);
    % \draw [semithick] (2.1,-.55) to [out=-90,in=90] (1,-1);
    % \draw (1,-1) node[below]{$k$};
    % \draw (1,1.4) node{};
    \end{tikzpicture}
    }\quad =\quad
    \hackcenter{
    \begin{tikzpicture}
    \draw [thick,->] (0,-.5) to [out=90,in=-90] (.5,.5);
    \draw [thick,->] (.5,-.5)to [out=90,in=-90] (1,.5);
    \draw [thick,->] (2,-.5) to [out=90,in=-90] (0,.5);
    \draw [thick,->] (1.5,-.5) to [out=90,in=-90] (2,.5);
    % \fill [white] (-.25,-.25) rectangle (2.25,.25);
    % \draw [thick] (-.25,-.25) rectangle (2.25,.25);
    % \draw (1,0) node{$k$};
    \draw (1.15,-.3) node{$\cdots$};
%    \draw (1.25,.4) node{$\cdots$};
    \draw [semithick] (-.1,-.55) to [out=-90,in=90] (1,-1);
    \draw [semithick] (2.1,-.55) to [out=-90,in=90] (1,-1);
    \draw (1,-1) node[below]{$k$};
    % \draw (1,-1.2) node{};
    \draw (1,1.4) node{};
    \end{tikzpicture}
    }\, .
\end{equation*}
In the case when a partition $\pi =(\pi_{1},\pi_{2},\dots,\pi_{l})\in \cP$ is given, we align the corresponding boxes to draw
\begin{equation*}
    \hackcenter{
    \begin{tikzpicture}
    \draw [thick,->] (0,-.5)--(0,.5);
    \draw [thick,->] (.5,-.5)--(.5,.5);
    \draw [thick,->] (2,-.5)--(2,.5);
    \fill [lightgray] (-.25,-.25) rectangle (2.25,.25);
    \draw [thick] (-.25,-.25) rectangle (2.25,.25);
    \draw (1,0) node{$\pi$};
    \draw (1.25,-.4) node{$\cdots$};
    \draw (1.25,.4) node{$\cdots$};
    % \draw [semithick] (-.1,-.55) to [out=-90,in=90] (1,-1);
    % \draw [semithick] (2.1,-.55) to [out=-90,in=90] (1,-1);
%    \draw (1,-1) node[below]{$k$};
%    \draw (1,1.4) node{};
    \end{tikzpicture}
    }\quad =\quad
    \hackcenter{
    \begin{tikzpicture}
    \draw [thick,->] (0,-.5)--(0,.5);
    \draw [thick,->] (1,-.5)--(1,.5);
    \fill [white] (-.2,-.25) rectangle (1.2,.25);
    \draw [thick] (-.2,-.25) rectangle (1.2,.25);
    \draw (.5,0) node{$\pi_{1}$};
    \draw (.5,-.4) node{$\cdots$};
    \draw (.5,.4) node{$\cdots$};
    \draw [thick,->] (1.5,-.5)--(1.5,.5);
    \draw [thick,->] (2.5,-.5)--(2.5,.5);
    \fill [white] (1.3,-.25) rectangle (2.7,.25);
    \draw [thick] (1.3,-.25) rectangle (2.7,.25);
    \draw (2,0) node{$\pi_{2}$};
    \draw (2,-.4) node{$\cdots$};
    \draw (2,.4) node{$\cdots$};
    \draw (3.25,0) node{$\cdots$};
    \draw [thick,->] (4,-.5)--(4,.5);
    \draw [thick,->] (5,-.5)--(5,.5);
    \fill [white] (3.8,-.25) rectangle (5.2,.25);
    \draw [thick] (3.8,-.25) rectangle (5.2,.25);
    \draw (4.5,0) node{$\pi_{l}$};
    \draw (4.5,-.4) node{$\cdots$};
    \draw (4.5,.4) node{$\cdots$};
    \end{tikzpicture}
    }\, .
\end{equation*}
Finally, again for a partition $\pi\in \cP$, we set
\begin{equation*}
    \alpha_{\pi}=
    \hackcenter{
    \begin{tikzpicture}
    % \draw [thick] (-1.5,-.5)--(-1.5,.5);
    \draw [thick,->] (-1.5,.25) to [out=90,in=90] (1.5,0);
    \draw [thick] (1.5,0) to [out=-90,in=-90] (-1.5,-.25);
    % \draw [thick,->] (.5,-.5)--(.5,.5);
    % \draw [thick] (-.5,-.5)--(-.5,.5);
    \draw [thick,->] (-.5,.25) to [out=90,in=90] (.5,0);
    \draw [thick] (.5,0) to [out=-90,in=-90] (-.5,-.25);
    \fill [lightgray] (-1.75,-.25) rectangle (-.25,.25);
    \draw [thick] (-1.75,-.25) rectangle (-.25,.25);
    \draw (-1,0) node{$\pi$};
    \draw (-1,-.4) node{$\cdots$};
    \draw (-1,.4) node{$\cdots$};
    \draw (1,0) node{$\cdots$};
    \end{tikzpicture}
    }\, \in Z(\cH).
\end{equation*}

The following result from \cite{KLM2019} is what we rely on in the proof of our main result.
\begin{thm}[\cite{KLM2019}]
\label{thm:KLM-theorem}
Under the homomorphism $\cF \colon Z(\cH)\to \mathrm{Fun} (\bY)$,
we have
\begin{equation*}
    \cF (c_{k})=B_{k+2},\quad k\in\bZ_{\geq 0},
\end{equation*}
and
\begin{equation*}
    \cF (\alpha_{\pi})=\Sigma_{\pi},\quad \pi\in \cP .
\end{equation*}
\end{thm}

\begin{rem}
Together with Proposition~\ref{prop:center_generators} and the fact that Boolean cumulants are algebraically independent, Theorem~\ref{thm:KLM-theorem} implies that the homomorphism $\cF$ is injective.
\end{rem}

\begin{rem}
In \cite{KLM2019}, it was established that $Z(\cH)$ is isomorphic to the ring of shifted symmetric functions.
For our purpose, however, the embedding of $Z(\cH)$ into $\mathrm{Fun} (\bY)$ is sufficient.
\end{rem}

\begin{rem}
In \cite{KLM2019}, it was also shown that $\cF (\tilde{c}_{k})=M_{k}$, $k\in \bZ_{\geq 0}$.
In fact, the relations between $c_{k}$, $k\in \bZ_{\geq 0}$ and $\tilde{c}_{k}$, $k\in \bZ_{\geq 0}$ seen in Proposition~\ref{prop:center_generators} are equivalent to the functional relation
\begin{equation*}
    G_{\lambda}(z)H_{\lambda}(z)=1,\quad \lambda\in \bY
\end{equation*}
satisfied by the generating functions of the moments and the Boolean cumulants.
\end{rem}

\section{Proof of Theorems \texorpdfstring{\ref{thm:main1}}{Main1} and \texorpdfstring{\ref{thm:main2}}{Main2}}
\label{sect:proofs_main1_main2}
In this section, we prove Theorems \ref{thm:main1} and \ref{thm:main2} at the same time.

\subsection{Alternative statement}
Due to Theorem~\ref{thm:KLM-theorem},
the first half of our main result Theorem~\ref{thm:main1} is equivalent to that
for any $\pi\in \cP$, $(-1)^{\ell (\pi)}\alpha_{\pi}$ is expressed by a polynomial of $-c_{k}$, $k=0,1,\dots$ with non-negative integer coefficients.

Here we show a slightly stronger result than Theorem~\ref{thm:main1}.
For a partition $\pi\in \cP_{n}$ and a sequence $\bm{i}=(i_{1},i_{2},\dots,i_{n})\in (\bZ_{\geq 0})^{n}$, we set
\begin{equation*}
    \alpha_{\pi}(\bm{i})=
    \hackcenter{
    \begin{tikzpicture}
    \draw [thick] (-3,.25)--(-3,.5);
    \draw [thick,->] (-3,.5) to [out=90,in=90] (3,0);
    \draw [thick] (3,0) to [out=-90,in=-90] (-3,-.25);
    \fill [black] (-3,.5) circle (.1);
    \draw (-3,.5) node[left]{$i_{1}$};
    \draw [thick] (-2.25,.25)--(-2.25,.5);
    \draw [thick,->] (-2.25,.5) to [out=90,in=90] (2.25,0);
    \draw [thick] (2.25,0) to [out=-90,in=-90] (-2.25,-.25);
    \fill [black] (-2.25,.5) circle (.1);
    \draw (-2.25,.5) node[left]{$i_{2}$};
    \draw [thick] (-1,.25)--(-1,.5);
    \draw [thick,->] (-1,.5) to [out=90,in=90] (1,0);
    \draw [thick] (1,0) to [out=-90,in=-90] (-1,-.25);
    \fill [black] (-1,.5) circle (.1);
    \draw (-1,.5) node[left]{$i_{n}$};
    \fill [lightgray] (-3.25,-.25) rectangle (-.75,.25);
    \draw [thick] (-3.25,-.25) rectangle (-.75,.25);
    \draw (-2,0) node{$\pi$};
    \draw (-1.75,-.5) node{$\cdots$};
    \draw (-1.75,.5) node{$\cdots$};
    \draw (1.5,0) node{$\cdots$};
    \end{tikzpicture}
    }\quad .
\end{equation*}

We prove the following theorem in this section.
\begin{thm}
\label{thm:main_alt}
For each $\pi\in \cP$ and a sequence $\bm{i}\in (\bZ_{\geq 0})^{|\pi|}$,
there exists a polynomial $P_{\pi,\bm{i}}(y_{0},y_{1},\dots)$ with non-negative integer coefficients such that
\begin{equation*}
    (-1)^{\ell (\pi)}\alpha_{\pi,\bm{i}}=P_{\pi,\bm{i}}(-c_{0},-c_{1},\dots)
\end{equation*}
holds in $Z(\cH)$.
Furthermore, when we define degrees of variables as $\deg y_{i}=i$, $i=0,1,\dots$, all nonzero monomials appearing in $P_{\pi,\bm{i}}$ are of degrees at most $|\pi|-\ell (\pi)+|\bm{i}|$, where $|\bm{i}|=\sum_{k=1}^{|\pi|}i_{k}$.
\end{thm}

Specializing $\bm{i}=(0,\dots,0)$ in Theorem~\ref{thm:main_alt},
we prove Theorem~\ref{thm:main1}.

\subsection{Proof of Theorem~\ref{thm:main_alt}}
We prove it by induction in terms of $|\pi|$.
In the base case $|\pi|=1$, the only possibility of $\pi$ is $\pi=(1)$. Hence for any $\bm{i}=(i)$ with $i\in \bZ_{\geq 0}$,
\begin{equation*}
    \alpha_{(1)}(\bm{i})=
    \hackcenter{
    \begin{tikzpicture}
    \draw [thick,->] (-.5,0) arc (180:0:.5);
    \draw [thick] (.5,0) arc (360:180:.5);
    \fill [black] (-.5,0) circle (.1);
    \draw (-.5,0) node[left]{$i$};
    \end{tikzpicture}
    }\,=c_{i}.
\end{equation*}
Therefore, we can take $P_{(1),(i)}=y_{i}$ for each $i\in \bZ_{\geq 0}$,
which is certainly of degree $i$.

Let us assume that the assertion is true for all $\pi$ up to $|\pi|\leq n$,
and prove that the assertion is again true for $\pi\in \cP_{n+1}$.

We consider two cases separately.
The first case is that $\pi=(n+1)$ is of length $1$.
In this case, we observe for any sequence $\bm{i}=(i_{1},\dots,i_{n},i_{n+1})\in (\bZ_{\geq 0})^{n+1}$ that
\begin{equation*}
    \alpha_{\pi}(\bm{i})=
    \hackcenter{
    \begin{tikzpicture}
    \draw [thick] (-3.25,-.25) rectangle (-0.75,.25);
    \draw [thick] (-3,.25) -- (-3,.5);
    \draw [thick] (-1,.25) -- (-1,.5);
    \draw [thick] (-1.5,.25) -- (-1.5,.5);
    \draw [thick] (-1,-.75) to [out=90,in=-90] (-.5,-.25);
    \draw [thick] (-.5,-.75) to [out=90,in=-90] (-1,-.25);
    % \draw [thick,->] (-.5,-.25) arc (180:0:.5);
    \draw [thick,->] (-.5,-.25) to [out=90,in=90] (.5,-.25);
    \draw [thick] (.5,-.25) -- (.5,-.75);
    \draw [thick] (.5,-.75) to [out=-90,in=-90] (-.5,-.75);
    % \draw [thick] (.5,-.75) arc (360:180:.5);
    %
    \draw [thick,->] (-1,.5) to [out=90,in=90] (1,0);
    \draw [thick] (1,0)--(1,-.75);
    \draw [thick] (1,-.75) to [out=-90,in=-90] (-1,-.75);
    \draw [thick,->] (-3,.5) to [out=90,in=90] (3,0);
    \draw [thick] (3,0) to [out=-90,in=-90] (-3,-.25);
    \draw [thick,->] (-1.5,.5) to [out=90,in=90] (1.5,0);
    \draw [thick] (1.5,0) -- (1.5,-.75);
    \draw [thick] (1.5,-.75) to [out=-90,in=-90] (-1.5,-.75);
    \draw [thick] (-1.5,-.75)--(-1.5,-.25);
    \fill [black] (-3,.5) circle (.1);
    \fill [black] (-1.5,.5) circle (.1);
    \fill [black] (-1,.5) circle (.1);
    \fill [black] (-.5,-.25) circle (.1);
    \draw (-2,0) node{$n$};
    \draw (-3,.5) node[left]{$i_{1}$};
    \draw (-2.25,-.5) node{$\cdots$};
    \draw (-1.5,.5) node[left]{$i_{n-1}$};
    \draw (-1,.5) node[right]{$i_{n}$};
    \draw (-.5,.-.25) node[right]{$i_{n+1}$};
    \draw (2.25,0) node{$\cdots$};
    \end{tikzpicture}
    }\quad .
\end{equation*}
Applying Lemma~\ref{lem:dot_moving_over_crossing} repeatedly to move the dots located on the most internal curl over the crossing, we have
\begin{align*}
    (-1)\alpha_{\pi}(\bm{i})=&
    (-1)\hackcenter{
    \begin{tikzpicture}
    \draw [thick] (-3.25,-.25) rectangle (-0.75,.25);
    \draw [thick] (-3,.25) -- (-3,.5);
    \draw [thick] (-1,.25) -- (-1,.5);
    \draw [thick] (-1.5,.25) -- (-1.5,.5);
    \draw [thick,->] (-1,.5) to [out=90,in=90] (1,0);
    \draw [thick] (1,0) to [out=-90,in=-90] (-1,-.25);
    \draw [thick,->] (-3,.5) to [out=90,in=90] (3,0);
    \draw [thick] (3,0) to [out=-90,in=-90] (-3,-.25);
    \draw [thick,->] (-1.5,.5) to [out=90,in=90] (1.5,0);
    \draw [thick] (1.5,0) to [out=-90,in=-90] (-1.5,-.25);
    \fill [black] (-3,.5) circle (.1);
    \fill [black] (-1.5,.5) circle (.1);
    \fill [black] (-1,.5) circle (.1);
    \draw (-2,0) node{$n$};
    \draw (-3,.5) node[left]{$i_{1}$};
    \draw (-2.25,-.5) node{$\cdots$};
    \draw (-1.5,.5) node[left]{$i_{n-1}$};
    \draw [blue, thick, dashed,->] (0,2.5) ..controls (-1,2.5) and (0,.5) ..(-.85,.5);
    \draw (0,2.5) node[right]{$i_{n}+i_{n+1}+1$};
    \draw (2.25,0) node{$\cdots$};
    \end{tikzpicture}
    } \\
    &+\sum_{b=0}^{i_{n+1}-1}
    \hackcenter{
    \begin{tikzpicture}
    \draw [thick] (-3.25,-.25) rectangle (-0.75,.25);
    \draw [thick] (-3,.25) -- (-3,.5);
    \draw [thick] (-1,.25) -- (-1,.5);
    \draw [thick] (-1.5,.25) -- (-1.5,.5);
    \draw [thick,->] (-1,.5) to [out=90,in=90] (1,0);
    \draw [thick] (1,0) to [out=-90,in=-90] (-1,-.25);
    \draw [thick,->] (-.5,0) arc (180:0:.5);
    \draw [thick] (.5,0) arc (360:180:.5);
    \draw [thick,->] (-3,.5) to [out=90,in=90] (3,0);
    \draw [thick] (3,0) to [out=-90,in=-90] (-3,-.25);
    \draw [thick,->] (-1.5,.5) to [out=90,in=90] (1.5,0);
    \draw [thick] (1.5,0) to [out=-90,in=-90] (-1.5,-.25);
    \fill [black] (-3,.5) circle (.1);
    \fill [black] (-1.5,.5) circle (.1);
    \fill [black] (-1,.5) circle (.1);
    \fill [black] (-.5,0) circle (.1);
    \draw (-2,0) node{$n$};
    \draw (-3,.5) node[left]{$i_{1}$};
    \draw (-2.25,-.5) node{$\cdots$};
    \draw (-1.5,.5) node[left]{$i_{n-1}$};
    \draw [blue, thick, dashed,->] (0,2.5) ..controls (-1,2.5) and (0,.5) ..(-.85,.5);
    \draw (0,2.5) node[right]{$i_{n}+i_{n+1}-b-1$};
    \draw (-.5,0) node[right]{$b$};
    \draw (2.25,0) node{$\cdots$};
    \end{tikzpicture}
    }\quad .
\end{align*}
By the induction hypothesis, there exists a non-negative integer coefficient polynomial $\wtilde{P}_{\pi,\bm{i}}(y_{0},y_{1},\dots)$ such that the first term in the right hand side is expressed as
\begin{equation*}
    \wtilde{P}_{\pi,\bm{i}}(-c_{0},-c_{1},\dots).
\end{equation*}
Furthermore, each monomial appearing in $\wtilde{P}_{\pi,\bm{i}}$ is of degree at most $n+|\bm{i}|$.

Let us consider each summand in the second term.
We shall apply Proposition~\ref{prop:bubble_move_non-negativity} repeatedly to move the internal hoop outside.
For each $b=0,1,\dots, i_{n+1}-1$, we have
\begin{align*}
    &\hackcenter{
    \begin{tikzpicture}
    \draw [thick] (-3.25,-.25) rectangle (-0.75,.25);
    \draw [thick] (-3,.25) -- (-3,.5);
    \draw [thick] (-1,.25) -- (-1,.5);
    \draw [thick] (-1.5,.25) -- (-1.5,.5);
    \draw [thick,->] (-1,.5) to [out=90,in=90] (1,0);
    \draw [thick] (1,0) to [out=-90,in=-90] (-1,-.25);
    \draw [thick,->] (-.5,0) arc (180:0:.5);
    \draw [thick] (.5,0) arc (360:180:.5);
    \draw [thick,->] (-3,.5) to [out=90,in=90] (3,0);
    \draw [thick] (3,0) to [out=-90,in=-90] (-3,-.25);
    \draw [thick,->] (-1.5,.5) to [out=90,in=90] (1.5,0);
    \draw [thick] (1.5,0) to [out=-90,in=-90] (-1.5,-.25);
    \fill [black] (-3,.5) circle (.1);
    \fill [black] (-1.5,.5) circle (.1);
    \fill [black] (-1,.5) circle (.1);
    \fill [black] (-.5,0) circle (.1);
    \draw (-2,0) node{$n$};
    \draw (-3,.5) node[left]{$i_{1}$};
    \draw (-2.25,-.5) node{$\cdots$};
    \draw (-1.5,.5) node[left]{$i_{n-1}$};
    \draw [blue, thick, dashed,->] (0,2.5) ..controls (-1,2.5) and (0,.5) ..(-.85,.5);
    \draw (0,2.5) node[right]{$i_{n}+i_{n+1}-b-1$};
    \draw (-.5,0) node[right]{$b$};
    \draw (2.25,0) node{$\cdots$};
    \end{tikzpicture}
    }\\
    =&\sum_{\substack{\pr{\bm{i}}\in (\bZ_{\geq 0})^{n},j\in \bZ_{\geq 0}\\ |\pr{\bm{i}}|+j\leq |\bm{i}|-1}}m_{\pr{\bm{i}},j}^{b}
    \hackcenter{
    \begin{tikzpicture}
    \draw [thick] (-3.25,-.25) rectangle (-0.75,.25);
    \draw [thick] (-3,.25) -- (-3,.5);
    \draw [thick] (-1,.25) -- (-1,.5);
    \draw [thick] (-1.5,.25) -- (-1.5,.5);
    \draw [thick,->] (-1,.5) to [out=90,in=90] (1,0);
    \draw [thick] (1,0) to [out=-90,in=-90] (-1,-.25);
    \draw [thick,->] (4,0) arc (180:0:.5);
    \draw [thick] (5,0) arc (360:180:.5);
    \draw [thick,->] (-3,.5) to [out=90,in=90] (3,0);
    \draw [thick] (3,0) to [out=-90,in=-90] (-3,-.25);
    \draw [thick,->] (-1.5,.5) to [out=90,in=90] (1.5,0);
    \draw [thick] (1.5,0) to [out=-90,in=-90] (-1.5,-.25);
    \fill [black] (-3,.5) circle (.1);
    \fill [black] (-1.5,.5) circle (.1);
    \fill [black] (-1,.5) circle (.1);
    \fill [black] (4,0) circle (.1);
    \draw (-2,0) node{$n$};
    \draw (-3,.5) node[left]{$\pr{i}_{1}$};
    \draw (-2.25,-.5) node{$\cdots$};
    \draw (-1.5,.5) node[left]{$\pr{i}_{n-1}$};
    % \draw [blue, thick, dashed,->] (0,2.5) ..controls (-1,2.5) and (0,.5) ..(-.85,.5);
    \draw (-1,.5) node[right]{$\pr{i}_{n}$};
    \draw (4,0) node[left]{$j$};
    \draw (2.25,0) node{$\cdots$};
    \end{tikzpicture}
    } \\
    &-\sum_{\substack{\ppr{\bm{i}}\in (\bZ_{\geq 0})^{n}\\ |\ppr{\bm{i}}|\leq |\bm{i}|-1}}n_{\ppr{\bm{i}}}^{b}
    \hackcenter{
    \begin{tikzpicture}
    \draw [thick] (-3.25,-.25) rectangle (-0.75,.25);
    \draw [thick] (-3,.25) -- (-3,.5);
    \draw [thick] (-1,.25) -- (-1,.5);
    \draw [thick] (-1.5,.25) -- (-1.5,.5);
    \draw [thick,->] (-1,.5) to [out=90,in=90] (1,0);
    \draw [thick] (1,0) to [out=-90,in=-90] (-1,-.25);
    %
    % \draw [thick,->] (4,0) arc (180:0:.5);
    % \draw [thick] (5,0) arc (360:180:.5);
    %
    \draw [thick,->] (-3,.5) to [out=90,in=90] (3,0);
    \draw [thick] (3,0) to [out=-90,in=-90] (-3,-.25);
    \draw [thick,->] (-1.5,.5) to [out=90,in=90] (1.5,0);
    \draw [thick] (1.5,0) to [out=-90,in=-90] (-1.5,-.25);
    \fill [black] (-3,.5) circle (.1);
    \fill [black] (-1.5,.5) circle (.1);
    \fill [black] (-1,.5) circle (.1);
    % \fill [black] (4,0) circle (.1);
    %
    \draw (-2,0) node{$n$};
    \draw (-3,.5) node[left]{$\ppr{i}_{1}$};
    \draw (-2.25,-.5) node{$\cdots$};
    \draw (-1.5,.5) node[left]{$\ppr{i}_{n-1}$};
    % \draw [blue, thick, dashed,->] (0,2.5) ..controls (-1,2.5) and (0,.5) ..(-.85,.5);
    \draw (-1,.5) node[right]{$\ppr{i}_{n}$};
    % \draw (4,0) node[left]{$c$};
    \draw (2.25,0) node{$\cdots$};
    \end{tikzpicture}
    }\, ,
\end{align*}
where $m_{\pr{i},j}^{b}, n_{\ppr{\bm{i}}}^{b}\in\bZ_{\geq 0}$
for all $\pr{\bm{i}}, \ppr{\bm{i}}\in (\bZ_{\geq 0})^{n}$, $j\in \bZ_{\geq 0}$.
By the induction hypothesis, there exists a non-negative integer coefficient polynomial $\what{P}_{\pi,\bm{i}}^{b}(y_{0},y_{1},\dots)$
such that the above element is expressed as
\begin{equation*}
    \what{P}^{b}_{\pi,\bm{i}}(-c_{0},-c_{1},\dots).
\end{equation*}
Furthermore, all monomials appearing in $\what{P}^{b}_{\pi,\bm{i}}$ are of degrees at most $n+|\bm{i}|-2$.

Consequently, when we set
\begin{equation*}
    P_{\pi,\bm{i}}=\wtilde{P}_{\pi,\bm{i}}+\sum_{b=0}^{i_{n+1}-1}\what{P}^{b}_{\pi,\bm{i}},
\end{equation*}
it is a non-negative integer coefficient polynomial of $y_{0},y_{1},\dots$
consisting of monomials of degrees at most $n+|\bm{i}|$, such that
\begin{equation*}
    (-1)\alpha_{\pi}(\bm{i})=P_{\pi,\bm{i}}(-c_{0},-c_{1},\dots).
\end{equation*}

The second case is that $\ell=\ell (\pi)\geq 2$.
In this case, setting $\wtilde{\pi}=(\pi_{1},\dots,\pi_{\ell -1})$ and
$k=|\wtilde{\pi}|$, we can draw, for $\bm{i}\in (\bZ_{\geq 0})^{n+1}$
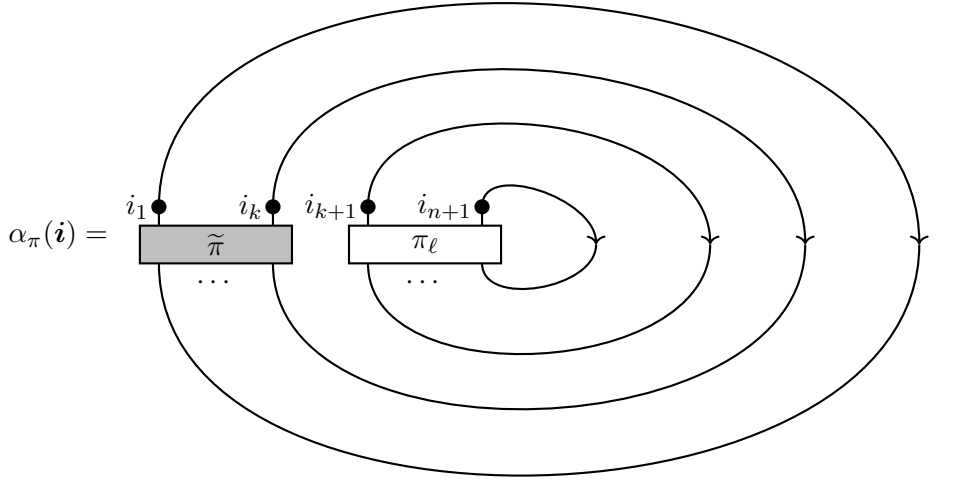
\begin{equation*}
    \alpha_{\pi}(\bm{i})=
    \hackcenter{
    \begin{tikzpicture}
    \draw [thick] (-5,.25)--(-5,.5);
    \draw [thick,->] (-5,.5) to [out=90,in=90] (5,0);
    \draw [thick] (5,0) to [out=-90,in=-90] (-5,-.25);
    \fill [black] (-5,.5) circle (.1);
    \draw (-5,.5) node[left]{$i_{1}$};
    \draw [thick] (-3.5,.25)--(-3.5,.5);
    \draw [thick,->] (-3.5,.5) to [out=90,in=90] (3.5,0);
    \draw [thick] (3.5,0) to [out=-90,in=-90] (-3.5,-.25);
    \fill [black] (-3.5,.5) circle (.1);
    \draw (-3.5,.5) node[left]{$i_{k}$};
    \draw [thick] (-2.25,.25)--(-2.25,.5);
    \draw [thick,->] (-2.25,.5) to [out=90,in=90] (2.25,0);
    \draw [thick] (2.25,0) to [out=-90,in=-90] (-2.25,-.25);
    \fill [black] (-2.25,.5) circle (.1);
    \draw (-2.25,.5) node[left]{$i_{k+1}$};
    \draw [thick] (-.75,.25)--(-.75,.5);
    \draw [thick,->] (-.75,.5) to [out=90,in=90] (.75,0);
    \draw [thick] (.75,0) to [out=-90,in=-90] (-.75,-.25);
    \fill [black] (-.75,.5) circle (.1);
    \draw (-.75,.5) node[left]{$i_{n+1}$};
    \fill [lightgray] (-5.25,-.25) rectangle (-3.25,.25);
    \draw [thick] (-5.25,-.25) rectangle (-3.25,.25);
    \draw (-4.25,0) node{$\wtilde{\pi}$};
    \draw (-4.25,-.5) node{$\cdots$};
    %
    % \fill [white] (-5.25,-.25) rectangle (-3.25,.25);
    \draw [thick] (-2.5,-.25) rectangle (-.5,.25);
    \draw (-1.5,0) node{$\pi_{\ell}$};
    \draw (-1.5,-.5) node{$\cdots$};
    % \draw (-1.75,.5) node{$\cdots$};
    % \draw (1.5,0) node{$\cdots$};
    \end{tikzpicture}
    }\quad .
\end{equation*}
Since $\pi_{\ell}\leq n$, we can apply the induction hypothesis to the most internal block. Specifically, there exists a non-negative integer coefficient polynomial $\pr{P}(y_{0},y_{1},\dots)$ such that
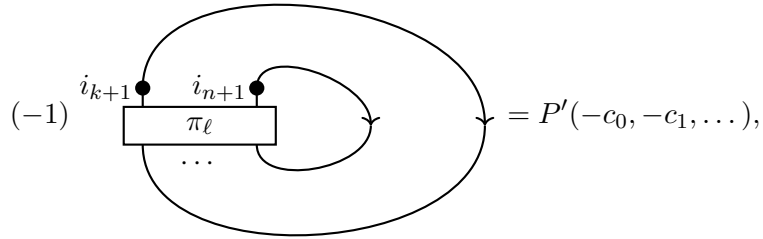
\begin{equation*}
    (-1)\hackcenter{
    \begin{tikzpicture}
    \draw [thick] (-2.25,.25)--(-2.25,.5);
    \draw [thick,->] (-2.25,.5) to [out=90,in=90] (2.25,0);
    \draw [thick] (2.25,0) to [out=-90,in=-90] (-2.25,-.25);
    \fill [black] (-2.25,.5) circle (.1);
    \draw (-2.25,.5) node[left]{$i_{k+1}$};
    \draw [thick] (-.75,.25)--(-.75,.5);
    \draw [thick,->] (-.75,.5) to [out=90,in=90] (.75,0);
    \draw [thick] (.75,0) to [out=-90,in=-90] (-.75,-.25);
    \fill [black] (-.75,.5) circle (.1);
    \draw (-.75,.5) node[left]{$i_{n+1}$};
    \draw [thick] (-2.5,-.25) rectangle (-.5,.25);
    \draw (-1.5,0) node{$\pi_{\ell}$};
    \draw (-1.5,-.5) node{$\cdots$};
    \end{tikzpicture}
    }\,= \pr{P}(-c_{0},-c_{1},\dots),
\end{equation*}
and all monomials appearing in $\pr{P}$ are of degrees at most $\pi_{\ell}+\sum_{s=k+1}^{n+1}i_{s}-1$.
% Here we need to recall the trace $\Tr (\cH)$ and its action on $Z(\cH)$.
% When we write
% \begin{equation*}
%     f_{\pi}(\bm{i})=(-1)^{\ell -1}
%     \hackcenter{
%     \begin{tikzpicture}
%     \draw [thick,->] (-.75,-.5)--(-.75,.75);
%     \draw [thick,->] (.75,-.5)--(.75,.75);
%     \fill [black] (-.75,.5) circle (.1);
%     \fill [black] (.75,.5) circle (.1);
%     \fill [lightgray] (-1,-.25) rectangle (1,.25);
%     \draw [thick] (-1,-.25) rectangle (1,.25);
%     \draw (0,0) node{$\widetilde{\pi}$};
%     \draw (-.75,.5) node[left]{$i_{1}$};
%     \draw (.75,.5) node[left]{$i_{k}$};
%     \end{tikzpicture}
%     }\quad \in \End_{\cH}(Q_{(k)}),
% \end{equation*}
% we can identify $(-1)^{\ell}\alpha_{\pi}(\bm{i})$ as
% \begin{equation*}
%     (-1)^{\ell}\alpha_{\pi}(\bm{i})=[f_{\pi}(\bm{i})].\pr{P}(-c_{0},-c_{1},\dots)
% \end{equation*}
% by means of the action of $\Tr (\cH)$ on $Z(\cH)$.
Therefore, $(-1)^{\ell}\alpha_{\pi}(\bm{i})$
is a linear combination with coefficients being non-negative integers of terms of the form
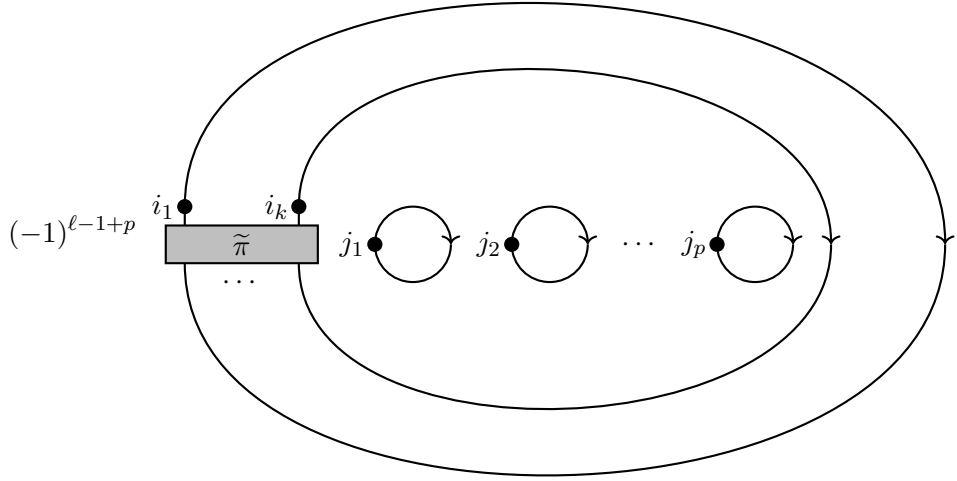
\begin{equation*}
    (-1)^{\ell-1+p}\,\hackcenter{
    \begin{tikzpicture}
    \draw [thick] (-5,.25)--(-5,.5);
    \draw [thick,->] (-5,.5) to [out=90,in=90] (5,0);
    \draw [thick] (5,0) to [out=-90,in=-90] (-5,-.25);
    \fill [black] (-5,.5) circle (.1);
    \draw (-5,.5) node[left]{$i_{1}$};
    \draw [thick] (-3.5,.25)--(-3.5,.5);
    \draw [thick,->] (-3.5,.5) to [out=90,in=90] (3.5,0);
    \draw [thick] (3.5,0) to [out=-90,in=-90] (-3.5,-.25);
    \fill [black] (-3.5,.5) circle (.1);
    \draw (-3.5,.5) node[left]{$i_{k}$};
   \draw [thick,->] (-2.5,0) arc (180:0:.5);
   \draw [thick] (-1.5,0) arc (360:180:.5);
   \fill [black] (-2.5,0) circle (.1);
   \draw (-2.5,0) node[left]{$j_{1}$};
   \draw [thick,->] (-.7,0) arc (180:0:.5);
   \draw [thick] (.3,0) arc (360:180:.5);
   \fill [black] (-.7,0) circle (.1);
   \draw (-.7,0) node[left]{$j_{2}$};
   \draw (1,0) node{$\cdots$};
   \draw [thick,->] (2,0) arc (180:0:.5);
   \draw [thick] (3,0) arc (360:180:.5);
   \fill [black] (2,0) circle (.1);
   \draw (2,0) node[left]{$j_{p}$};
    % % 
    \fill [lightgray] (-5.25,-.25) rectangle (-3.25,.25);
    \draw [thick] (-5.25,-.25) rectangle (-3.25,.25);
    \draw (-4.25,0) node{$\wtilde{\pi}$};
    \draw (-4.25,-.5) node{$\cdots$};
    \end{tikzpicture}
    }
\end{equation*}
for various $\bm{j}=(j_{1},\dots,j_{p})\in (\bZ_{\geq 0})^{p}$, $p=0,1,\dots$ such that $|\bm{j}|:=\sum_{r=1}^{p}j_{r}\leq \pi_{\ell}+\sum_{s=k+1}^{n+1}i_{s}-1$.

We can use Proposition~\ref{prop:bubble_move_non-negativity}
to move bubbles inside outside in the above expression for each $\bm{j}$ to see that $(-1)^{\ell}\alpha_{\pi}(\bm{i})$ becomes
\begin{align*}
    \sum_{q\geq 0}\sum_{\substack{\pr{\bm{i}}\in (\bZ_{\geq 0})^{k},\pr{\bm{j}}\in (\bZ_{\geq 0})^{q}\\ |\pr{\bm{i}}|+|\pr{\bm{j}}|\leq \pi_{\ell}+|\bm{i}|-1}}(-1)^{\ell-1+q}m_{\pr{\bm{i}},\pr{\bm{j}}}\quad
    \hackcenter{
    \begin{tikzpicture}
    \draw [thick] (-2.5,.25)--(-2.5,.5);
    \draw [thick,->] (-2.5,.5) to [out=90,in=90] (2,0);
    \draw [thick] (2,0) to [out=-90,in=-90] (-2.5,-.25);
    \fill [black] (-2.5,.5) circle (.1);
    \draw (-2.5,.5) node[left]{$\pr{i}_{1}$};
    \draw [thick] (-1,.25)--(-1,.5);
    \draw [thick,->] (-1,.5) to [out=90,in=90] (1,0);
    \draw [thick] (1,0) to [out=-90,in=-90] (-1,-.25);
    \fill [black] (-1,.5) circle (.1);
    \draw (-1,.5) node[left]{$\pr{i}_{k}$};
   \draw [thick,->] (3,1) arc (180:0:.5);
  \draw [thick] (4,1) arc (360:180:.5);
  \fill [black] (3,1) circle (.1);
  \draw (3,1) node[left]{$\pr{j}_{1}$};
  \draw [thick,->] (3,-1) arc (180:0:.5);
  \draw [thick] (4,-1) arc (360:180:.5);
  \fill [black] (3,-1) circle (.1);
  \draw (3,-1) node[left]{$\pr{j}_{q}$};
   \draw (3.5,0) node{$\vdots$};
   %
%   \draw [thick,->] (2,0) arc (180:0:.5);
%   \draw [thick] (3,0) arc (360:180:.5);
%   \fill [black] (2,0) circle (.1);
%   \draw (2,0) node[left]{$j_{p}$};
    % % 
    \fill [lightgray] (-2.75,-.25) rectangle (-.75,.25);
    \draw [thick] (-2.75,-.25) rectangle (-.75,.25);
    \draw (-1.75,0) node{$\wtilde{\pi}$};
    \draw (-1.75,-.5) node{$\cdots$};
    \end{tikzpicture}
    }\, ,
\end{align*}
where $m_{\pr{\bm{i}},\pr{\bm{j}}}\in \bZ_{\geq 0}$ for all $\pr{\bm{i}}, \pr{\bm{j}}$.
Since $|\wtilde{\pi}|\leq n$, we can apply the induction hypothesis;
for each $\pr{\bm{i}}\in (\bZ_{\geq 0})^{k}$,
there exists a non-negative coefficient polynomial $P_{\wtilde{\pi},\pr{\bm{i}}}(y_{0},y_{1},\dots)$ such that
\begin{equation*}
    (-1)^{\ell -1}\hackcenter{
    \begin{tikzpicture}
    \draw [thick] (-2.5,.25)--(-2.5,.5);
    \draw [thick,->] (-2.5,.5) to [out=90,in=90] (2,0);
    \draw [thick] (2,0) to [out=-90,in=-90] (-2.5,-.25);
    \fill [black] (-2.5,.5) circle (.1);
    \draw (-2.5,.5) node[left]{$\pr{i}_{1}$};
    \draw [thick] (-1,.25)--(-1,.5);
    \draw [thick,->] (-1,.5) to [out=90,in=90] (1,0);
    \draw [thick] (1,0) to [out=-90,in=-90] (-1,-.25);
    \fill [black] (-1,.5) circle (.1);
    \draw (-1,.5) node[left]{$\pr{i}_{k}$};
    \fill [lightgray] (-2.75,-.25) rectangle (-.75,.25);
    \draw [thick] (-2.75,-.25) rectangle (-.75,.25);
    \draw (-1.75,0) node{$\wtilde{\pi}$};
    \draw (-1.75,-.5) node{$\cdots$};
    \end{tikzpicture}
    }\, =P_{\wtilde{\pi},\pr{\bm{i}}}(-c_{0},-c_{1},\dots ),
\end{equation*}
and all monomials are of degrees at most $k-(\ell -1)+|\pr{\bm{i}}|$.

Let us set
\begin{equation*}
    P_{\pi,\bm{i}}(y_{0},y_{1},\dots)=\sum_{q\geq 0}\sum_{\substack{\pr{\bm{i}}\in (\bZ_{\geq 0})^{k},\pr{\bm{j}}\in (\bZ_{\geq 0})^{q}\\ |\pr{\bm{i}}|+|\pr{\bm{j}}|\leq \pi_{\ell}+|\bm{i}|-1}}m_{\pr{\bm{i}},\pr{\bm{j}}}P_{\wtilde{\pi},\pr{\bm{i}}}(y_{0},y_{1},\dots)y_{\pr{j}_{1}}\cdots y_{\pr{j}_{q}}.
\end{equation*}
By construction, all monomials appearing in $P_{\pi,\bm{i}}$ are of degrees at most
\begin{equation*}
    k-(\ell-1)+\pi_{\ell}+|\bm{i}|-1=|\pi|-\ell (\pi)+|\bm{i}|,
\end{equation*}
and the desired identity
\begin{equation*}
    (-1)^{\ell}\alpha_{\pi}(\bm{i})=P_{\pi,\bm{i}}(-c_{0},-c_{1},\dots )
\end{equation*}
holds.

Therefore, the assertion is true for $\pi\in \cP_{n+1}$, completing the proof.

\section{Proof of Theorem \texorpdfstring{\ref{thm:main3}}{Main3}}
\label{sect:proof_main3}
\subsection{Preliminary arguments}
If $\sigma\in S_{n}$, we take the longest element, i.e., the inversion
\begin{equation*}
    w_{0}=\left(
    \begin{array}{cccc}
    1 & 2 &\cdots & n \\
    n & n-1 &\cdots & 1
    \end{array}
    \right)\in S_{n}
\end{equation*}
in the same group, and set
\begin{equation*}
    \sigma^{*}:=w_{0}\sigma^{-1}w_{0}.
\end{equation*}
It is clear that $(\sigma^{*})^{*}=\sigma$.
Now, if we recall the notational convention introduced after Proposition~\ref{prop:symmetric_group_algebra_Hom_space}, we may observe, for instance, the following: for $\sigma,\rho\in S_{n}$,
\begin{align*}
\hackcenter{
\begin{tikzpicture}
    \draw[thick,->] (-2.25,.25) to [out=90,in=180] (0,1.5);
    \draw[thick] (0,1.5) to [out=0, in=90] (2.25,.25);
    \draw[thick,->] (2.25,-.25) to [out=-90,in=0] (0,-1.5);
    \draw[thick] (0,-1.5) to [out=180,in=-90] (-2.25,-.25);
    \draw[thick,->] (-.75,.25) to [out=90,in=180] (0,.75);
    \draw[thick] (0,.75) to [out=0, in=90] (.75,.25);
    \draw[thick,->] (.75,-.25) to [out=-90,in=0] (0,-.75);
    \draw[thick] (0,-.75) to [out=180,in=-90] (-.75,-.25);
    \fill[gray] (-2.5,-.25) rectangle (-.5,.25);
    \draw[thick] (-2.5,-.25) rectangle (-.5,.25);
    \draw (-1.5,0) node{$\sigma$};
    \fill[gray] (.5,-.25) rectangle (2.5,.25);
    \draw[thick] (.5,-.25) rectangle (2.5,.25);
    \draw (1.5,0) node{$\rho$};
    \draw ($(-2.25,-.4)!.5!(-.75,-.4)$) node{$\cdots$};
    \draw ($(2.25,-.4)!.5!(.75,-.4)$) node{$\cdots$};
    \end{tikzpicture}
    }\quad &=\quad
    \hackcenter{
    \begin{tikzpicture}
    \draw[thick] (-2.25,-.1) to (-2.25,.1);
    \draw[thick] (-.75,-.1) to (-.75,.1);
    \draw ($(-2.25,-.8)!.5!(-.75,-.8)$) node{$\cdots$};
    \draw[thick,->] (-2.25,.6) to [out=90,in=90] (2.25,0);
    \draw[thick] (2.25,0) to [out=-90,in=-90] (-2.25,-.6);
    \draw[thick,->] (-.75,.6) to [out=90,in=90] (.75,0);
    \draw[thick] (.75,0) to [out=-90,in=-90] (-.75,-.6);
    \draw ($(2.25,0)!.5!(.75,0)$) node{$\cdots$};
    \fill[gray] (-2.5,-.6) rectangle (-.5,-.1);
    \draw[thick] (-2.5,-.6) rectangle (-.5,-.1);
    \draw ($(-2.5,-.6)!.5!(-.5,-.1)$) node{$\sigma$};
    \fill[gray] (-2.5,.1) rectangle (-.5,.6);
    \draw[thick] (-2.5,.1) rectangle (-.5,.6);
    \draw ($(-2.5,.1)!.5!(-.5,.6)$) node{$\rho^{*}$};
\end{tikzpicture}
    }\\
    &=\hackcenter{
    \begin{tikzpicture}
    \draw[thick] (-2.25,-.1) to (-2.25,.1);
    \draw[thick] (-.75,-.1) to (-.75,.1);
    \draw ($(-2.25,-.8)!.5!(-.75,-.8)$) node{$\cdots$};
    \draw[thick,->] (-2.25,.6) to [out=90,in=90] (2.25,0);
    \draw[thick] (2.25,0) to [out=-90,in=-90] (-2.25,-.6);
    \draw[thick,->] (-.75,.6) to [out=90,in=90] (.75,0);
    \draw[thick] (.75,0) to [out=-90,in=-90] (-.75,-.6);
    \draw ($(2.25,0)!.5!(.75,0)$) node{$\cdots$};
    \fill[gray] (-2.5,-.6) rectangle (-.5,-.1);
    \draw[thick] (-2.5,-.6) rectangle (-.5,-.1);
    \draw ($(-2.5,-.6)!.5!(-.5,-.1)$) node{$\rho^{*}$};
    \fill[gray] (-2.5,.1) rectangle (-.5,.6);
    \draw[thick] (-2.5,.1) rectangle (-.5,.6);
    \draw ($(-2.5,.1)!.5!(-.5,.6)$) node{$\sigma$};
\end{tikzpicture}
    }
\end{align*}
by continuously deforming the diagram, under which taking $*$ is realized by the 180-degree rotation.

Let us suppose that $\sigma\in S_{n}$ belongs to the conjugacy class corresponding to $\pi\in \cP_{n}$.
Then, there exists $\rho\in S_{n}$ such that
\begin{equation*}
    \rho\sigma\rho^{-1}=(1\dots\pi_{1})(\pi_{1}+1\dots \pi_{1}+\pi_{2})\cdots (\pi_{1}+\cdots \pi_{\ell(\pi)-1}\dots n).
\end{equation*}
The above observation allows us to get 
\begin{align}
\label{eq:norm_char_diagram_conjugate}
    \hackcenter{
    \begin{tikzpicture}
    \draw[thick,->] (-2.25,.25) to [out=90,in=90] (2.25,0);
    \draw[thick] (2.25,0) to [out=-90,in=-90] (-2.25,-.25);
    \draw[thick,->] (-.75,.25) to [out=90,in=90] (.75,0);
    \draw[thick] (.75,0) to [out=-90,in=-90] (-.75,-.25);
    \fill[gray] (-2.5,-.25) rectangle (-.5,.25);
    \draw[thick] (-2.5,-.25) rectangle (-.5,.25);
    \draw (-1.5,0) node{$\sigma$};
    \draw ($(-2,-.4)!.5!(-.75,-.4)$) node{$\cdots$};
    \end{tikzpicture}
    }\quad =\quad
    \hackcenter{
    \begin{tikzpicture}
    \draw[thick,->] (-2.25,.25) to [out=90,in=180] (0,1);
    \draw[thick] (0,1) to [out=0,in=90] (2.25,.6);
    \draw[thick,->] (2.25,-.6) to [out=-90,in=0] (0,-1);
    \draw[thick] (0,-1) to [out=180,in=-90] (-2.25,-.25);
    \draw[thick,->] (-.75,.25) to [out=90,in=180] (0,.8);
    \draw[thick] (0,.8) to [out=0,in=90] (.75,.6);
    \draw[thick,->] (.75,-.6) to [out=-90,in=0] (0,-.8);
    \draw[thick] (0,-.8) to [out=180,in=-90] (-.75,-.25);
    \draw[thick] (.75,.6)--(.75,-.6);
    \draw[thick] (2.25,.6)--(2.25,-.6);
    \fill[lightgray] (-2.5,-.25) rectangle (-.5,.25);
    \draw[thick] (-2.5,-.25) rectangle (-.5,.25);
    \draw ($(-2.5,-.25)!.5!(-.5,.25)$) node{$\pi$};
    \draw ($(-2,-.4)!.5!(-.75,-.4)$) node{$\cdots$};
    \fill[gray] (.5,.1) rectangle (2.5,.6);
    \draw[thick] (.5,.1) rectangle (2.5,.6);
    \fill[gray] (.5,-.1) rectangle (2.5,-.6);
    \draw[thick] (.5,-.1) rectangle (2.5,-.6);
    \draw ($(.5,-.75)!.5!(2.5,-.75)$)node{$\cdots$};
    \draw ($(.5,.1)!.5!(2.5,.6)$)node{$\rho^{*}$};
    \draw ($(.5,-.1)!.5!(2.5,-.6)$)node{$(\rho^{-1})^{*}$};
    \end{tikzpicture}
    }\quad =\alpha_{\pi}
\end{align}
as $\rho^{*}(\rho^{-1})^{*}=e$.

In the following statement, if $\sigma\in S_{n}$, we agree that $|\sigma|=n$.
\begin{prop}
\label{prop:colliding_permutations}
Let $\sigma\in S_{n}$ and $\rho\in S_{m}$. 
We also take $a\leq n$ and $b\leq m$.
Then, there exists an expansion
\begin{align*}
    &\hackcenter{
    \begin{tikzpicture}
    \draw[thick,->] (-3.9,-1)--(-3.9,1);
    \draw[thick,->] (-3.2,-1)--(-3.2,1);
    \draw ($(-3.9,-.5)!.5!(-3.2,-.5)$)node{$\cdots$};
    \draw [semithick] (-4,-1.1) to [out=-90,in=90] (-3.55,-1.25);
    \draw [semithick] (-3.1,-1.1) to [out=-90,in=90] (-3.55,-1.25);
    \draw (-3.55,-1.25) node[below]{$a$};
    \fill[gray] (-4,-.25) rectangle (-2,.25);
    \draw[thick] (-4,-.25) rectangle (-2,.25);
    \draw (-3,0) node{$\sigma$};
    \draw[thick,->] (-3,.25) to [out=90,in=90] (-.2,0);
    \draw[thick] (-.2,0) to [out=-90,in=-90] (-3,-.25);
    \draw[thick,->] (-2.1,.25) to [out=90,in=90] (-.9,0);
    \draw[thick] (-.9,0) to [out=-90,in=-90] (-2.1,-.25);
    \draw ($(-3,-.5)!.5!(-2.1,-.5)$)node{$\cdots$};
    \draw[thick,<-] (3.9,-1)--(3.9,1);
    \draw[thick,<-] (3.2,-1)--(3.2,1);
    \draw ($(3.9,-.5)!.5!(3.2,-.5)$)node{$\cdots$};
    \draw [semithick] (4,-1.1) to [out=-90,in=90] (3.55,-1.25);
    \draw [semithick] (3.1,-1.1) to [out=-90,in=90] (3.55,-1.25);
    \draw (3.55,-1.25) node[below]{$b$};
    \fill[gray] (4,-.25) rectangle (2,.25);
    \draw[thick] (4,-.25) rectangle (2,.25);
    \draw (3,0) node{$\rho$};
    \draw[thick] (3,.25) to [out=90,in=90] (.2,0);
    \draw[thick,<-] (.2,0) to [out=-90,in=-90] (3,-.25);
    \draw[thick] (2.1,.25) to [out=90,in=90] (.9,0);
    \draw[thick,<-] (.9,0) to [out=-90,in=-90] (2.1,-.25);
    \draw ($(3,-.5)!.5!(2.1,-.5)$)node{$\cdots$};
    \end{tikzpicture}
    } \\
    =&\sum_{\sigma',\rho'}\bfm_{\sigma',\rho'}\quad
    \hackcenter{
    \begin{tikzpicture}
    \draw[thick,->] (-2.4,-1)--(-2.4,1);
    \draw[thick,->] (-1.6,-1)--(-1.6,1);
    \draw ($(-2.4,-.5)!.5!(-1.6,-.5)$)node{$\cdots$};
    \fill[gray] (-2.5,-.25) rectangle (-.5,.25);
    \draw[thick] (-2.5,-.25) rectangle (-.5,.25);
    \draw ($(-2.5,0)!.5!(-.5,0)$) node{$\sigma'$};
    \draw [semithick] (-2.5,-1.1) to [out=-90,in=90] (-2,-1.25);
    \draw [semithick] (-1.5,-1.1) to [out=-90,in=90] (-2,-1.25);
    \draw (-2,-1.25) node[below]{$a$};
    \draw[thick,<-] (2.4,-1)--(2.4,1);
    \draw[thick,<-] (1.6,-1)--(1.6,1);
    \draw ($(2.4,-.5)!.5!(1.6,-.5)$)node{$\cdots$};
    \fill[gray] (2.5,-.25) rectangle (.5,.25);
    \draw[thick] (2.5,-.25) rectangle (.5,.25);
    \draw ($(2.5,0)!.5!(.5,0)$) node{$\rho'$};
    \draw [semithick] (2.5,-1.1) to [out=-90,in=90] (2,-1.25);
    \draw [semithick] (1.5,-1.1) to [out=-90,in=90] (2,-1.25);
    \draw (2,-1.25) node[below]{$b$};
    \draw[thick,->] (-1.4,.25) to [out=90,in=180] (0,1);
    \draw[thick] (0,1) to [out=0,in=90] (1.4,.25);
    \draw[thick,->] (1.4,-.25) to [out=-90,in=0] (0,-1);
    \draw[thick] (0,-1) to [out=180,in=-90] (-1.4,-.25);
    \draw[thick,->] (-.6,.25) to [out=90,in=180] (0,.5);
    \draw[thick] (0,.5) to [out=0,in=90] (.6,.25);
    \draw[thick,->] (.6,-.25) to [out=-90,in=0] (0,-.5);
    \draw[thick] (0,-.5) to [out=180,in=-90] (-.6,-.25);
    \draw ($(-1.4,-.5)!.5!(-.6,-.5)$) node{$\cdots$};
    \end{tikzpicture}
    }\quad,
\end{align*}
where the sum runs over permutations $\sigma'$ and $\rho'$ such that $|\sigma'|-a=|\rho'|-b$ and $|\rho'|\leq |\sigma|+|\rho|-a$  (equivalently, $|\sigma'|\leq |\sigma|+|\rho|-b$ under the first constraint) and all coefficients $\bfm_{\sigma',\rho'}$ are non-negative integers.
Furthermore, $\bfm_{\sigma',\rho'}=0$ unless $\sgn (\sigma')=\sgn (\sigma)$ and $\sgn (\rho')=\sgn (\rho)$ are satisfied.
\end{prop}
\begin{proof}
Let us first comment on the range of summation in the right-hand side.
The condition $|\sigma'|-a=|\rho'|-b$ is required simply to make sense of the diagram.
The condition $|\rho'| \leq |\sigma|+|\rho|-a$ means that the number of strings entering $\rho'$ does not exceed the number of strings originally entering $\rho$ and those leaving and returning to $|\sigma|$.

The proof goes by induction on $|\sigma|-a$.
If $a=|\sigma|$, the left-hand side is already in the desired form by setting $\rho'=\rho$ and $\sigma'$ to be the image of $\sigma$ along the embedding $S_{n}\hookrightarrow S_{n+|\rho|-b}$ that fixes the last $|\rho|-b$ letters. In this case, $\sgn (\sigma')=\sgn (\sigma)$ and $\sgn (\rho')=\sgn (\rho)$ are indeed satisfied.

Assume that the claim is true when $a\geq k+1$ with some $k< |\sigma|$.
Then, by repeating the application of the second relation of (\ref{eq:hom_relations_res_ind}), we get
\begin{align*}
    &\hackcenter{
    \begin{tikzpicture}
    \draw[thick,->] (-3.9,-1)--(-3.9,1);
    \draw[thick,->] (-3.2,-1)--(-3.2,1);
    \draw ($(-3.9,-.5)!.5!(-3.2,-.5)$)node{$\cdots$};
    \draw [semithick] (-4,-1.1) to [out=-90,in=90] (-3.55,-1.25);
    \draw [semithick] (-3.1,-1.1) to [out=-90,in=90] (-3.55,-1.25);
    \draw (-3.55,-1.25) node[below]{$k$};
    \fill[gray] (-4,-.25) rectangle (-2,.25);
    \draw[thick] (-4,-.25) rectangle (-2,.25);
    \draw (-3,0) node{$\sigma$};
    \draw[thick,->] (-3,.25) to [out=90,in=90] (-.2,0);
    \draw[thick] (-.2,0) to [out=-90,in=-90] (-3,-.25);
    \draw[thick,->] (-2.1,.25) to [out=90,in=90] (-.9,0);
    \draw[thick] (-.9,0) to [out=-90,in=-90] (-2.1,-.25);
    \draw[thick,->] (-2.8,.25) to [out=90,in=90] (-.4,0);
    \draw[thick] (-.4,0) to [out=-90,in=-90] (-2.8,-.25);
    \draw ($(-2.8,-.5)!.5!(-2.1,-.5)$)node{$\cdots$};
    \draw[thick,<-] (3.9,-1)--(3.9,1);
    \draw[thick,<-] (3.2,-1)--(3.2,1);
    \draw ($(3.9,-.5)!.5!(3.2,-.5)$)node{$\cdots$};
    \draw [semithick] (4,-1.1) to [out=-90,in=90] (3.55,-1.25);
    \draw [semithick] (3.1,-1.1) to [out=-90,in=90] (3.55,-1.25);
    \draw (3.55,-1.25) node[below]{$b$};
    \fill[gray] (4,-.25) rectangle (2,.25);
    \draw[thick] (4,-.25) rectangle (2,.25);
    \draw (3,0) node{$\rho$};
    \draw[thick] (3,.25) to [out=90,in=90] (.2,0);
    \draw[thick,<-] (.2,0) to [out=-90,in=-90] (3,-.25);
    \draw[thick] (2.1,.25) to [out=90,in=90] (.9,0);
    \draw[thick,<-] (.9,0) to [out=-90,in=-90] (2.1,-.25);
    \draw ($(3,-.5)!.5!(2.1,-.5)$)node{$\cdots$};
    \end{tikzpicture}
    } \\
    =\quad &\hackcenter{
    \begin{tikzpicture}
    \draw[thick,->] (-3.9,-1)--(-3.9,1);
    \draw[thick,->] (-3.2,-1)--(-3.2,1);
    \draw[thick] (-3,-1)--(-3,1);
    \draw ($(-3.9,-.5)!.5!(-3.2,-.5)$)node{$\cdots$};
    \draw [semithick] (-4,-1.1) to [out=-90,in=90] (-3.55,-1.25);
    \draw [semithick] (-3.1,-1.1) to [out=-90,in=90] (-3.55,-1.25);
    \draw (-3.55,-1.25) node[below]{$k$};
    \fill[gray] (-4,-.25) rectangle (-2,.25);
    \draw[thick] (-4,-.25) rectangle (-2,.25);
    \draw (-3,0) node{$\sigma$};
    \draw[thick,->] (-2.1,.25) to [out=90,in=90] (-.9,0);
    \draw[thick] (-.9,0) to [out=-90,in=-90] (-2.1,-.25);
    \draw[thick,->] (-2.8,.25) to [out=90,in=90] (-.4,0);
    \draw[thick] (-.4,0) to [out=-90,in=-90] (-2.8,-.25);
    \draw ($(-2.8,-.5)!.5!(-2.1,-.5)$)node{$\cdots$};
    \draw[thick,->] (-3,1) to [out=90,in=180] (-2.5,1.2);
    \draw[thick,->] (-2.5,1.2)--(2.5,1.2);
    \draw[thick] (2.5,1.2) to [out=0,in=90] (3,1)--(3,.7) to [out=-90, in=90] (1.9,.25) -- (1.9,-.25) to [out=-90,in=90] (3,-.7)--(3,-1);
    \draw[thick,->] (3,-1) to [out=-90,in=0] (2.5,-1.2);
    \draw[thick,->] (2.5,-1.2)--(-2.5,-1.2);
    \draw[thick] (-2.5,-1.2) to [out=180,in=-90] (-3,-1);
    \draw[thick,<-] (3.9,-1)--(3.9,1);
    \draw[thick,<-] (3.2,-1)--(3.2,1);
    \draw ($(3.9,-.5)!.5!(3.2,-.5)$)node{$\cdots$};
    \draw [semithick] (4,-1.1) to [out=-90,in=90] (3.55,-1.25);
    \draw [semithick] (3.1,-1.1) to [out=-90,in=90] (3.55,-1.25);
    \draw (3.55,-1.25) node[below]{$b$};
    \fill[gray] (4,-.25) rectangle (2,.25);
    \draw[thick] (4,-.25) rectangle (2,.25);
    \draw (3,0) node{$\rho$};
    \draw[thick] (3,.25) to [out=90,in=90] (.2,0);
    \draw[thick,<-] (.2,0) to [out=-90,in=-90] (3,-.25);
    \draw[thick] (2.1,.25) to [out=90,in=90] (.9,0);
    \draw[thick,<-] (.9,0) to [out=-90,in=-90] (2.1,-.25);
    \draw ($(3,-.6)!.5!(2.1,-.6)$)node{$\cdots$};
    \draw[thick, blue, dashed] (-4.2,1.1)--(4.2,1.1)--(4.2,-1.1)--(-4.2,-1.1)--cycle;
    \end{tikzpicture}
    } \\
    & + \sum_{i=1}^{|\rho|-b}\quad 
    \hackcenter{
    \begin{tikzpicture}
    \draw[thick,->] (-3.9,-1)--(-3.9,1);
    \draw[thick,->] (-3.2,-1)--(-3.2,1);
    \draw[thick] (-3,-1)--(-3,1);
    \draw ($(-3.9,-.5)!.5!(-3.2,-.5)$)node{$\cdots$};
    \draw [semithick] (-4,-1.1) to [out=-90,in=90] (-3.55,-1.25);
    \draw [semithick] (-3.1,-1.1) to [out=-90,in=90] (-3.55,-1.25);
    \draw (-3.55,-1.25) node[below]{$k$};
    \fill[gray] (-4,-.25) rectangle (-2,.25);
    \draw[thick] (-4,-.25) rectangle (-2,.25);
    \draw (-3,0) node{$\sigma$};
    \draw[thick,->] (-2.1,.25) to [out=90,in=90] (-.9,0);
    \draw[thick] (-.9,0) to [out=-90,in=-90] (-2.1,-.25);
    \draw[thick,->] (-2.8,.25) to [out=90,in=90] (-.4,0);
    \draw[thick] (-.4,0) to [out=-90,in=-90] (-2.8,-.25);
    \draw ($(-2.8,-.5)!.5!(-2.1,-.5)$)node{$\cdots$};
    \draw[thick,->] (-3,1) to [out=90,in=180] (-2.5,1.2);
    \draw[thick,->] (-2.5,1.2)--(2.5,1.2);
    \draw[thick] (2.5,1.2) to [out=0,in=90] (3,1)--(3,.7) to [out=-90, in=90] (2.5,.25) -- (2.5,-.25) to [out=-90,in=90] (3,-.7)--(3,-1);
    \draw[thick,->] (3,-1) to [out=-90,in=0] (2.5,-1.2);
    \draw[thick,->] (2.5,-1.2)--(-2.5,-1.2);
    \draw[thick] (-2.5,-1.2) to [out=180,in=-90] (-3,-1);
    \draw (2.5,.25)node[above]{\small $i$};
    \draw[thick,<-] (3.9,-1)--(3.9,1);
    \draw[thick,<-] (3.2,-1)--(3.2,1);
    \draw ($(3.9,-.5)!.5!(3.2,-.5)$)node{$\cdots$};
    \draw [semithick] (4,-1.1) to [out=-90,in=90] (3.55,-1.25);
    \draw [semithick] (3.1,-1.1) to [out=-90,in=90] (3.55,-1.25);
    \draw (3.55,-1.25) node[below]{$b$};
    \fill[gray] (4,-.25) rectangle (2,.25);
    \draw[thick] (4,-.25) rectangle (2,.25);
    \draw (3,0) node{$\rho$};
    \draw[thick] (3,.25) to [out=90,in=90] (.2,0);
    \draw[thick,<-] (.2,0) to [out=-90,in=-90] (3,-.25);
    \draw[thick] (2.1,.25) to [out=90,in=90] (.9,0);
    \draw[thick,<-] (.9,0) to [out=-90,in=-90] (2.1,-.25);
    \draw ($(3,-.6)!.5!(2.1,-.6)$)node{$\cdots$};
    \draw[thick, blue, dashed] (-4.2,1.1)--(4.2,1.1)--(4.2,-1.1)--(-4.2,-1.1)--cycle;
    \end{tikzpicture}
    }.
\end{align*}
Here we labeled the strings connected to $\rho$ by $i$ from left to right.
Here, we can apply the induction hypothesis locally to the blue dashed boxes in both the first term and the second sum.
If we keep track of the configuration of the strings, we can also see that the range of permutations appearing is as required.
Also, the permutations expressed by
\begin{equation*}
    \hackcenter{
    \begin{tikzpicture}
    %
%    \draw[thick,->] (-3,1) to [out=90,in=180] (-2.5,1.2);
%    \draw[thick,->] (-2.5,1.2)--(2.5,1.2);
    \draw[thick,->] (3,1)--(3,.7) to [out=-90, in=90] (1.9,.25) -- (1.9,-.25) to [out=-90,in=90] (3,-.7)--(3,-1);
%    \draw[thick,->] (3,-1) to [out=-90,in=0] (2.5,-1.2);
%    \draw[thick,->] (2.5,-1.2)--(-2.5,-1.2);
%    \draw[thick] (-2.5,-1.2) to [out=180,in=-90] (-3,-1);
    %
    \draw[thick,<-] (3.9,-1)--(3.9,1);
    \draw[thick,<-] (3.2,-1)--(3.2,1);
    \draw ($(3.9,-.5)!.5!(3.2,-.5)$)node{$\cdots$};
    \draw [semithick] (4,-1.1) to [out=-90,in=90] (3.55,-1.25);
    \draw [semithick] (3.1,-1.1) to [out=-90,in=90] (3.55,-1.25);
    \draw (3.55,-1.25) node[below]{$b$};
    \fill[gray] (4,-.25) rectangle (2,.25);
    \draw[thick] (4,-.25) rectangle (2,.25);
    \draw (3,0) node{$\rho$};
    \draw[thick] (3,.25) to [out=90,in=-90] (2.8,1);
    \draw[thick,<-] (2.8,-1) to [out=90,in=-90] (3,-.25);
    \draw[thick] (2.1,.25) to [out=90,in=-90] (1.9,1);
    \draw[thick,<-] (1.9,-1) to [out=90,in=-90] (2.1,-.25);
    \draw ($(3,-.6)!.5!(2.1,-.6)$)node{$\cdots$};
    \end{tikzpicture}
    }\qquad \mbox{and}\qquad
    \hackcenter{
    \begin{tikzpicture}
%    \draw[thick,->] (-3,1) to [out=90,in=180] (-2.5,1.2);
%    \draw[thick,->] (-2.5,1.2)--(2.5,1.2);
    \draw[thick,->] (3,1)--(3,.7) to [out=-90, in=90] (2.5,.25) -- (2.5,-.25) to [out=-90,in=90] (3,-.7)--(3,-1);
%    \draw[thick,->] (3,-1) to [out=-90,in=0] (2.5,-1.2);
%    \draw[thick,->] (2.5,-1.2)--(-2.5,-1.2);
%    \draw[thick] (-2.5,-1.2) to [out=180,in=-90] (-3,-1);
    \draw (2.5,.25)node[above]{\small $i$};
    \draw[thick,<-] (3.9,-1)--(3.9,1);
    \draw[thick,<-] (3.2,-1)--(3.2,1);
    \draw ($(3.9,-.5)!.5!(3.2,-.5)$)node{$\cdots$};
    \draw [semithick] (4,-1.1) to [out=-90,in=90] (3.55,-1.25);
    \draw [semithick] (3.1,-1.1) to [out=-90,in=90] (3.55,-1.25);
    \draw (3.55,-1.25) node[below]{$b$};
    \fill[gray] (4,-.25) rectangle (2,.25);
    \draw[thick] (4,-.25) rectangle (2,.25);
    \draw (3,0) node{$\rho$};
    \draw[thick] (3,.25) to [out=90,in=-90] (2.8,1);
    \draw[thick,<-] (2.8,-1) to [out=90,in=-90] (3,-.25);
    \draw[thick] (2.1,.25) to [out=90,in=-90] (2.1,1);
    \draw[thick,<-] (2.1,-1) to [out=90,in=-90] (2.1,-.25);
    \draw ($(3,-.6)!.5!(2.1,-.6)$)node{$\cdots$};
    \end{tikzpicture}
    }\quad, \quad i=1,\dots, |\rho|-b
\end{equation*}
have the same sign as $\rho$.
Therefore, the claim is true when $a=k$.
\end{proof}

If we apply Proposition~\ref{prop:colliding_permutations} to the case that $a=b=0$ and combine it with Theorem~\ref{thm:KLM-theorem} and (\ref{eq:norm_char_diagram_conjugate}), we get the following corollary (note also that if $\sigma$ belongs to the conjugacy class $\pi$, $\sgn (\sigma)=|\pi|-\ell (\pi)$):
\begin{cor}
\label{cor:product_norm_chars}
For $\pi_{1},\pi_{2}\in\cP$, the product $\Sigma_{\pi_{1}}\cdot \Sigma_{\pi_{2}}$ of the corresponding normalized characters is expanded in normalized characters $\Sigma_{\pi}$ such that $|\pi|\leq |\pi_{1}|+|\pi_{2}|$ and $(-1)^{|\pi|-\ell(\pi)}=(-1)^{|\pi_{1}|-\ell (\pi_{1})}(-1)^{|\pi_{2}|-\ell (\pi_{2})}$ with non-negative integer coefficients.
\end{cor}
This result in Corollary~\ref{cor:product_norm_chars} has been obtained in~\cite{IvanovKerov2001} with a combinatorial interpretation of the expansion coefficients.
It is natural to ask if our approach here recovers their combinatorial interpretation, but that would be beyond the scope of the present work.

\subsection{Alternative statement of Theorem~\ref{thm:main3}}
The strategy to prove Theorem~\ref{thm:main3} is similar to the proofs of Theorems~\ref{thm:main1} and \ref{thm:main2}.
That is, we employ the graphical calculus of the center of Khovanov's Heisenberg category.

Here we present a stronger result.
%\red{For $\sigma\in S_{n}$, we write $\ell (\sigma)$ for the Coxeter length of $\sigma$, i.e., the minimal number of simple transpositions $s_{i}$, $i=1,\dots, n-1$ needed to express $\sigma$.}
%Also recall the notational convention introduced after Proposition~\ref{prop:symmetric_group_algebra_Hom_space}.

\begin{thm}
\label{thm:main3_alt}
For each $k\in \bZ_{\geq 0}$, there exists an expansion
\begin{equation*}
    \hackcenter{
    \begin{tikzpicture}
    \draw[thick,->] (0,-1)--(0,1);
    \fill[black] (0,0) circle(.1);
    \draw (0,0)node[left]{$k$};
    \end{tikzpicture}
    }\, =
    \sum_{n=1}^{k+1}\sum_{\sigma\in S_{n}}\sfm^{k}_{\sigma}\quad
    \hackcenter{
    \begin{tikzpicture}
    \draw[thick,->] (-2.25,-1)--(-2.25,1);
    \draw[thick,->] (-.75,.25) to [out=90,in=90] (.75,0);
    \draw[thick] (.75,0) to [out=-90,in=-90] (-.75,-.25);
    \draw[thick,->] (-2,.25) to [out=90,in=90] (2,0);
    \draw[thick] (2,0) to [out=-90,in=-90] (-2,-.25);
    \fill[gray] (-2.5,-.25) rectangle (-.5,.25);
    \draw[thick] (-2.5,-.25) rectangle (-.5,.25);
    \draw (-1.5,0) node{$\sigma$};
    \draw ($(-2,-.4)!.5!(-.75,-.4)$) node{$\cdots$};
    \end{tikzpicture}
    }\, ,
\end{equation*}
where $\sfm^{k}_{\sigma}\in \bZ_{\geq 0}$ for all $\sigma$.
Furthermore $\sfm^{k}_{\sigma}=0$ unless $\sgn (\sigma) \equiv k$ mod $2$.
\end{thm}

Note that the expansion in Theorem~\ref{thm:main3_alt} is not unique.
In fact, for any $\sigma\in S_{n}$ and $\rho\in S_{n-1}$, we have the identity
\begin{equation*}
    \hackcenter{\begin{tikzpicture}
    \draw[thick,->] (-2.25,-1)--(-2.25,1);
    \draw[thick,->] (-.75,.25) to [out=90,in=90] (.75,0);
    \draw[thick] (.75,0) to [out=-90,in=-90] (-.75,-.25);
    \draw[thick,->] (-2,.25) to [out=90,in=90] (2,0);
    \draw[thick] (2,0) to [out=-90,in=-90] (-2,-.25);
    \fill[gray] (-2.5,-.25) rectangle (-.5,.25);
    \draw[thick] (-2.5,-.25) rectangle (-.5,.25);
    \draw (-1.5,0) node{$\sigma$};
    \draw ($(-2,-.4)!.5!(-.75,-.4)$) node{$\cdots$};
    \end{tikzpicture}}\quad =\quad
    \hackcenter{\begin{tikzpicture}
    \draw[thick,->] (-2.25,-1)--(-2.25,1);
    \draw[thick,->] (-.75,.8) to [out=90,in=90] (.75,0);
    \draw[thick] (.75,0) to [out=-90,in=-90] (-.75,-.8);
    \draw[thick] (-.75,-.8)--(-.75,.8);
    \draw[thick,->] (-2,.8) to [out=90,in=90] (2,0);
    \draw[thick] (2,0) to [out=-90,in=-90] (-2,-.8);
    \draw[thick] (-2,-.8)--(-2,.8);
    \fill[gray] (-2.5,-.25) rectangle (-.5,.25);
    \draw[thick] (-2.5,-.25) rectangle (-.5,.25);
    \fill[gray] (-2.1,.3) rectangle (-.5,.8);
    \draw[thick] (-2.1,.3) rectangle (-.5,.8);
    \fill[gray] (-2.1,-.8) rectangle (-.5,-.3);
    \draw[thick] (-2.1,-.8) rectangle (-.5,-.3);
    \draw (-1.5,0) node{$\sigma$};
    \draw ($(-2.1,.3)!.5!(-.5,.8)$) node{$\rho$};
    \draw ($(-2.1,-.8)!.5!(-.5,-.3)$) node{$\rho^{-1}$};
    \draw ($(-2,-1)!.5!(-.75,-1)$) node{$\cdots$};
    \end{tikzpicture}}\, .
\end{equation*}
Nevertheless, the fact that the coefficients are non-negative integers is independent of the way of expansion. Note also that the sign of a permutation is invariant under conjugation.

It is easy to deduce Theorem~\ref{thm:main3} from Theorem~\ref{thm:main3_alt} combined with Corollary~\ref{cor:product_norm_chars}.
\begin{proof}[Proof of Theorem~\ref{thm:main3}]
Owing to Corollary~\ref{cor:product_norm_chars}, it suffices to show that each $B_{k+2}$, $k\geq 0$ is expanded into $\Sigma_{\pi}$ such that $|\pi|\leq k+1$ with non-negative integer coefficients:
\begin{equation}
\label{eq:single_row_Boolean_expansion}
    B_{k+2}=\sum_{\pi\in\cP,\, |\pi|\leq k+1}m^{(k+2)}_{\pi}\Sigma_{\pi},\quad m^{(k+2)}_{\pi}\in\bZ_{\geq 0},
\end{equation}
where $m^{(k+2)}_{\pi}=0$ unless $|\pi|-\ell (\pi)\equiv k$ mod $2$.

For each $k\in\bZ$, we apply Theorem~\ref{thm:main3_alt} to the blue dashed box in
\begin{equation*}
    c_{k}=
    \hackcenter{
    \begin{tikzpicture}
    \draw[thick] (-.5,-.5) -- (-.5,.5);
    \draw[thick,->] (.5,.5) -- (.5,0);
    \draw[thick] (.5,0)--(.5,-.5);
    \draw[thick] (-.5,.5) arc (180:0:.5);
    \draw[thick] (.5,-.5) arc (360:180:.5);
    \fill[black] (-.5,0) circle(.1);
    \draw (-.5,0)node[left]{$k$};
    \draw[thick,blue,dashed] (-1,.5) -- (0,.5)--(0,-.5)--(-1,-.5)--cycle;
    \end{tikzpicture}
    }
\end{equation*}
to get an expansion
\begin{equation*}
    c_{k}=
    \sum_{n=1}^{k+1}\sum_{\sigma\in S_{n}}\sfm^{k}_{\sigma}\quad
    \hackcenter{
    \begin{tikzpicture}
    \draw[thick,->] (-2.25,.25) to [out=90,in=90] (2.25,0);
    \draw[thick] (2.25,0) to [out=-90,in=-90] (-2.25,-.25);
    \draw[thick,->] (-.75,.25) to [out=90,in=90] (.75,0);
    \draw[thick] (.75,0) to [out=-90,in=-90] (-.75,-.25);
    \fill[gray] (-2.5,-.25) rectangle (-.5,.25);
    \draw[thick] (-2.5,-.25) rectangle (-.5,.25);
    \draw (-1.5,0) node{$\sigma$};
    \draw ($(-2,-.4)!.5!(-.75,-.4)$) node{$\cdots$};
    \end{tikzpicture}
    }\, ,
\end{equation*}
where $\sfm^{k}_{\sigma}\in \bZ_{\geq 0}$ for all $\sigma$, and $\sfm^{k}_{\sigma}=0$ unless $\sgn (\sigma)\equiv k$ mod $2$.

% Suppose that $\sigma\in S_{n}$ belongs to the conjugacy class corresponding to $\pi\in \cP_{n}$.
% Then there exists $\rho\in S_{n}$ such that
% \begin{equation*}
%     \rho\sigma\rho^{-1}=(1\dots\pi_{1})(\pi_{1}+1\dots \pi_{1}+\pi_{2})\cdots (\pi_{1}+\cdots \pi_{\ell(\pi)-1}\dots n),
% \end{equation*}
% implying that
% \begin{equation*}
%     \hackcenter{
%     \begin{tikzpicture}
%     \draw[thick,->] (-2.25,.25) to [out=90,in=90] (2.25,0);
%     \draw[thick] (2.25,0) to [out=-90,in=-90] (-2.25,-.25);
%     \draw[thick,->] (-.75,.25) to [out=90,in=90] (.75,0);
%     \draw[thick] (.75,0) to [out=-90,in=-90] (-.75,-.25);
%     \fill[gray] (-2.5,-.25) rectangle (-.5,.25);
%     \draw[thick] (-2.5,-.25) rectangle (-.5,.25);
%     \draw (-1.5,0) node{$\sigma$};
%     \draw ($(-2,-.4)!.5!(-.75,-.4)$) node{$\cdots$};
%     \end{tikzpicture}
%     }\quad =\quad
%     \hackcenter{
%     \begin{tikzpicture}
%     \draw[thick,->] (-2.25,.25) to [out=90,in=90] (2.25,0);
%     \draw[thick] (2.25,0) to [out=-90,in=-90] (-2.25,-.25);
%     \draw[thick,->] (-.75,.25) to [out=90,in=90] (.75,0);
%     \draw[thick] (.75,0) to [out=-90,in=-90] (-.75,-.25);
%     \fill[gray] (-2.5,-.25) rectangle (-.5,.25);
%     \draw[thick] (-2.5,-.25) rectangle (-.5,.25);
%     \draw (-1.5,0) node{$\rho\sigma\rho^{-1}$};
%     \draw ($(-2,-.4)!.5!(-.75,-.4)$) node{$\cdots$};
%     \end{tikzpicture}
%     }\quad =\alpha_{\pi}.
% \end{equation*}
% Also note that in this case the Coxeter length of $\sigma$ is expressed as $\ell(\sigma)=|\pi|-\ell(\pi)$.
Now, (\ref{eq:norm_char_diagram_conjugate}) tells us that we can group the sum over $\sigma$ into the one over conjugacy classes.
For each $\pi\in \cP_{n}$, we set
\begin{equation*}
    m^{(k+2)}_{\pi}=\sum_{\substack{\sigma\in S_{n}\\\sigma\sim \pi}}\sfm^{k}_{\sigma}\in \bZ_{\geq 0},
\end{equation*}
where we wrote $\sigma\sim \pi$ to express that $\sigma$ belongs to the conjugacy class $\pi$.
As we have already mentioned, the sign of $\sigma$ is expressed as $|\pi|-\ell (\pi)$ using its conjugacy class $\pi$. Therefore, we have $m^{(k+2)}_{\pi}=0$ unless $|\pi|-\ell (\pi)\equiv k$ mod $2$.
Then, under the correspondence in Theorem~\ref{thm:KLM-theorem},
we arrive at the desired expansion (\ref{eq:single_row_Boolean_expansion}).
\end{proof}

\subsection{Proof of Theorem \texorpdfstring{\ref{thm:main3_alt}}{Main3alt}}
We prove Theorem~\ref{thm:main3_alt} by induction in $k$.
When $k=0$, the assertion is obviously true.

Assume that the assertion is true for all $k$ up to $n$.
Applying Lemma~\ref{lem:dot_moving_over_crossing} repeatedly, we have
\begin{equation}
\label{eq:dotted_arrow_induction}
    \hackcenter{
    \begin{tikzpicture}
    \draw [thick,->] (0,-1)--(0,1);
    \fill (0,0) circle(.1);
    \draw (0,0) node[left]{$n+1$};
    \end{tikzpicture}
    }\quad =\quad
    \hackcenter{
    \begin{tikzpicture}
    \draw [thick] (0,0) .. controls (.5,.7) and (1,.5) .. (1,0);
    \draw [thick] (0,0) .. controls (.5,-.7) and (1,-.5) .. (1,0);
    \draw [thick] (-.3,-1) .. controls (-.3,-1) and (-.3,-.4) ..(0,0);
    \draw [thick,->] (0,0) .. controls (-.3,.4) and (-.3,1) ..(-.3,1);
    \fill [black] (-.22,.5) circle (.1);
    \draw (-.22,.5) node[left]{$n$};
    \end{tikzpicture}
    }\quad =\quad
    \hackcenter{
    \begin{tikzpicture}
    \draw [thick] (0,0) .. controls (.5,.7) and (1,.5) .. (1,0);
    \draw [thick] (0,0) .. controls (.5,-.7) and (1,-.5) .. (1,0);
    \draw [thick] (-.3,-1) .. controls (-.3,-1) and (-.3,-.4) ..(0,0);
    \draw [thick,->] (0,0) .. controls (-.3,.4) and (-.3,1) ..(-.3,1);
    \fill [black] (.5,-.42) circle (.1);
    \draw (.5,-.42) node[below]{$n$};
    \end{tikzpicture}
    }\quad +\sum_{b=0}^{n-1}\quad
    \hackcenter{
    \begin{tikzpicture}
    \draw[thick,->] (-.5,-1)--(-.5,1);
    \fill (-.5,0) circle(.1);
    \draw (-.5,0) node[left]{$b$};
    \draw [thick,->] (0,0) arc (180:0:.5);
    \draw [thick] (1,0) arc (360:180:.5);
    \fill (.5,-.5) circle(.1);
    \draw (.5,-.5) node[below]{$n-1-b$};
    \end{tikzpicture}
    }
\end{equation}
Let us first study the first term.
By applying the induction hypothesis to the blue dashed box below, the first term becomes
\begin{equation*}
    \hackcenter{
    \begin{tikzpicture}
    \draw [thick] (-.5,-1)--(-.5,.5);
    \draw [thick] (-.5,.5) to [out=90,in=-90] (0,1);
    \draw [thick] (0,-.5)--(0,.5);
    \draw [thick,->] (0,.5) to [out=90, in=-90] (-.5,1)--(-.5,1.2);
    \draw[thick] (0,1) arc (180:0:.5);
    \draw[thick] (1,-.5) arc (360:180:.5);
    \draw[thick,->] (1,1)--(1,0);
    \draw[thick] (1,0)--(1,-.5);
    \fill [black] (0,0) circle (.1);
    \draw (0,0) node[right]{$n$};
    \draw[thick,blue,dashed] (-.25,.5)--(.5,.5)--(.5,-.5)--(-.25,-.5)--cycle;
    \end{tikzpicture}
    }\quad =
    \sum_{q=1}^{n+1}\sum_{\sigma\in S_{q}}\sfm^{n}_{\sigma}\quad
    \hackcenter{
    \begin{tikzpicture}
    \fill [gray] (-2.5,-.25) rectangle (-.5,.25);
    \draw [thick] (-2.5,-.25) rectangle (-.5,.25);
    \draw [thick,->] (-.75,.25) to [out=90,in=90] (.75,0);
    \draw [thick] (.75,0) to [out=-90,in=-90] (-.75,-.25);
    \draw [thick,->] (-2,.25) to [out=90,in=90] (2,0);
    \draw [thick] (2,0) to [out=-90,in=-90] (-2,-.25);
    \draw [thick] (-2.75,-1)--(-2.75,.25);
    \draw [thick] (-2.75,.25) to [out=90,in=-90] (-2.25,.75);
    \draw [thick] (-2.25,.25) to [out=90,in=-90] (-2.75,.75);
    \draw [thick,->] (-2.75,.75)--(-2.75,1);
    \draw [thick,->] (-2.25,.75) to [out=90,in=90] (2.25,0);
    \draw [thick] (2.25,0) to [out=-90,in=-90] (-2.25,-.25);
    \draw ($(-2.5,-.25)!.5!(-.5,.25)$) node{$\sigma$};
    \draw ($(-2,-.4)!.5!(-.75,-.4)$) node{$\cdots$};
    \end{tikzpicture}
    }\, ,
\end{equation*}
where $\sfm^{n}_{\sigma}\in \bZ_{\geq 0}$ and $\sfm^{n}_{\sigma}=0$ unless $\sgn (\sigma)\equiv n$ mod $2$.
Notice that this is already the desired form of expansion: if we think of $\sigma\in S_{q}$ as an element of $S_{q+1}$ permuting $\{2,\dots,q+1\}$ and compose the simple transpose $(1\,2)$ with it, we get a new permutation in $S_{q+1}$ whose sign is congruent to $n+1$ mod $2$. The above summation over $\sigma$ can be regarded as the summation over such obtained permutations.

Similarly, we can apply the induction hypothesis to each summand of the second term in~(\ref{eq:dotted_arrow_induction}).
\begin{equation*}
    \hackcenter{
    \begin{tikzpicture}
    \draw[thick,->] (-.5,-1)--(-.5,1);
    \fill (-.5,0) circle(.1);
    \draw (-.5,0) node[left]{$b$};
    \draw [thick,->] (0,0) arc (180:0:.5);
    \draw [thick] (1,0) arc (360:180:.5);
    \fill (.5,-.5) circle(.1);
    \draw (.5,-.5) node[below]{$n-1-b$};
    \end{tikzpicture}
    }=
    \sum_{p=1}^{b+1}\sum_{q=1}^{n-b}\sum_{\sigma\in S_{p}}\sum_{\rho\in S_{q}} \sfm^{b}_{\sigma} \sfm^{n-1-b}_{\rho}\quad
    \hackcenter{
    \begin{tikzpicture}
    \draw[thick,->] (-1.4,-1)--(-1.4,1);
    \draw[thick,->] (-.6,.25) to [out=90,in=90] (.6,0);
    \draw [thick] (.6,0) to [out=-90,in=-90] (-.6,-.25);
    \draw [thick,->] (-1.2,.25) to [out=90,in=90] (1.2,0);
    \draw [thick] (1.2,0) to [out=-90,in=-90] (-1.2,-.25);
    \fill[gray] (-1.5,-.25) rectangle (-.5,.25);
    \draw[thick] (-1.5,-.25) rectangle (-.5,.25);
    \draw ($(-1.5,-.25)!.5!(-.5,.25)$) node{$\sigma$};
    \draw ($(-1.2,-.4)!.5!(-.6,-.4)$) node{$\cdots$};
    \end{tikzpicture}
    }\quad 
    \hackcenter{
    \begin{tikzpicture}
    \draw[thick,->] (-.6,.25) to [out=90,in=90] (.6,0);
    \draw [thick] (.6,0) to [out=-90,in=-90] (-.6,-.25);
    \draw [thick,->] (-1.4,.25) to [out=90,in=90] (1.4,0);
    \draw [thick] (1.4,0) to [out=-90,in=-90] (-1.4,-.25);
    \fill[gray] (-1.5,-.25) rectangle (-.5,.25);
    \draw[thick] (-1.5,-.25) rectangle (-.5,.25);
    \draw ($(-1.5,-.25)!.5!(-.5,.25)$) node{$\rho$};
    \draw ($(-1.4,-.4)!.5!(-.6,-.4)$) node{$\cdots$};
    \end{tikzpicture}
    }\, ,
\end{equation*}
where $\sfm^{b}_{\sigma},\sfm^{n-1-b}_{\rho}\in \bZ_{\geq 0}$,
$\sfm^{b}_{\sigma}=0$ unless $\sgn (\sigma)\equiv b$ mod $2$, and $\sfm^{n-1-b}_{\rho}=0$ unless $\sgn (\rho)\equiv n-1-b$ mod $2$.
After the deformation
\begin{equation*}
    \hackcenter{
    \begin{tikzpicture}
    \draw[thick,->] (-.6,.25) to [out=90,in=90] (.6,0);
    \draw [thick] (.6,0) to [out=-90,in=-90] (-.6,-.25);
    \draw [thick,->] (-1.4,.25) to [out=90,in=90] (1.4,0);
    \draw [thick] (1.4,0) to [out=-90,in=-90] (-1.4,-.25);
    \fill[gray] (-1.5,-.25) rectangle (-.5,.25);
    \draw[thick] (-1.5,-.25) rectangle (-.5,.25);
    \draw ($(-1.5,-.25)!.5!(-.5,.25)$) node{$\rho$};
    \draw ($(-1.4,-.4)!.5!(-.6,-.4)$) node{$\cdots$};
    \end{tikzpicture}
    }\quad =\quad
    \hackcenter{
    \begin{tikzpicture}
    \draw[thick] (.6,.25) to [out=90,in=90] (-.6,0);
    \draw [thick,<-] (-.6,0) to [out=-90,in=-90] (.6,-.25);
    \draw [thick] (1.4,.25) to [out=90,in=90] (-1.4,0);
    \draw [thick,<-] (-1.4,0) to [out=-90,in=-90] (1.4,-.25);
    \fill[gray] (1.5,-.25) rectangle (.5,.25);
    \draw[thick] (1.5,-.25) rectangle (.5,.25);
    \draw ($(1.5,-.25)!.5!(.5,.25)$) node{$\rho^{*}$};
    \draw ($(1.4,-.4)!.5!(.6,-.4)$) node{$\cdots$};
    \end{tikzpicture}
    }\quad, 
\end{equation*}
we can apply Proposition~\ref{prop:colliding_permutations} in the case of $a=1$, $b=0$.
Consequently, we have the following expansion
\begin{equation*}
    \hackcenter{
    \begin{tikzpicture}
    \draw[thick,->] (-.5,-1)--(-.5,1);
    \fill (-.5,0) circle(.1);
    \draw (-.5,0) node[left]{$b$};
    \draw [thick,->] (0,0) arc (180:0:.5);
    \draw [thick] (1,0) arc (360:180:.5);
    \fill (.5,-.5) circle(.1);
    \draw (.5,-.5) node[below]{$n-1-b$};
    \end{tikzpicture}
    }=
    \sum_{q=1}^{n+1}\sum_{\sigma\in S_{q}}\wtilde{\sfm}^{b,n-1-b}_{\sigma}\quad
    \hackcenter{
    \begin{tikzpicture}
    \draw[thick,->] (-2.25,-1)--(-2.25,1);
    \draw[thick,->] (-.75,.25) to [out=90,in=90] (.75,0);
    \draw [thick] (.75,0) to [out=-90,in=-90] (-.75,-.25);
    \draw [thick,->] (-2,.25) to [out=90,in=90] (2,0);
    \draw [thick] (2,0) to [out=-90,in=-90] (-2,-.25);
    \fill[gray] (-2.5,-.25) rectangle (-.5,.25);
    \draw[thick] (-2.5,-.25) rectangle (-.5,.25);
    \draw ($(-2.5,-.25)!.5!(-.5,.25)$) node{$\sigma$};
    \draw ($(-2,-.4)!.5!(-.75,-.4)$) node{$\cdots$};
    \end{tikzpicture}
    }\, ,
\end{equation*}
where $\wtilde{\sfm}^{b,n-1-b}_{\sigma}\in\bZ_{\geq 0}$ and $\wtilde{\sfm}^{b,n-1-b}_{\sigma}=0$ unless $\sgn (\sigma)\equiv n-1\equiv n+1$ mod $2$.

Summing up all contributions, we obtain the desired expansion.

\section{Discussions}
\label{sect:discussion}
\subsection{Free cumulants}
In asymptotic representation theory,
instead of Boolean cumulants, free cumulants~\cite{Speicher1994} are more often discussed.
In contrast to that the generating function of Boolean cumulants is given by the multiplicative inverse of the Cauchy transform,
the generating function of free cumulants is defined by the inverse mapping of the Cauchy transform;
since the transition measure of a Young diagram $\lambda$ is compactly supported, its Cauchy transform $G_{\lambda}$ is a conformal map from a neighborhood of infinity to a neighborhood of the origin.
Its inverse map $K_{\lambda}$, often expressed as $K_{\lambda}=G_{\lambda}^{\langle -1\rangle}$ to distinguish it from the multiplicative inverse $H_{\lambda}=G_{\lambda}^{-1}$, is Laurent expanded around the origin as
\begin{equation*}
    K_{\lambda}(z)=\frac{1}{z}+\sum_{k=2}^{\infty}R_{k}(\lambda)z^{k-1}.
\end{equation*}
Here appearing $R_{k}(\lambda)$, $k=2,3\dots$ are the free cumulants of $\lambda$.

One importance of free cumulants in asymptotic representation theory lies in the fact that they describe the asymptotic behavior of characters of symmetric groups~\cite{Biane1998}.
Furthermore, it is known that the normalized characters corresponding to single-row partitions are expressed as polynomials (Kerov polynomials) of free cumulants $R_{k}$, $k=2,3,\dots$ with coefficients being non-negative integers~\cite{Biane2003,Feray2009,DolegaFeraySniady2010}.

Despite such importance of free cumulants, at the moment, graphically concise expressions of them are lacking.
It could be interesting to search for graphical expressions of free cumulants

\subsection{Combinatorial interpretation}
We have proved Theorems~\ref{thm:main1} and \ref{thm:main3} by induction in terms of the size of permutations or the number of dots, but we do not rely on any concrete combinatorics.
Nevertheless, the positivity of the coefficients strongly suggests the existence of underlying combinatorics.
Since, in principle, we can inductively compute expansion coefficients, one might hope that we could derive recursive formulas for those coeffieicnts, which will be helpful for combinatorial interpretation.
However, our induction bypasses diagrams that are neither normalized characters nor Boolean cumulants, which makes it hard to control the induction and get the recursive formulas.

Nevertheless, it would be interesting to give combinatorial proofs of  Theorems~\ref{thm:main1} and \ref{thm:main3}.
Note that the coefficients of Kerov polynomials have a clear combinatorial interpretation~\cite{Feray2009,DolegaFeraySniady2010}.

\subsection{Other variants, deformations}
There are quantum deformation~\cite{LicataSavage2013} and higher level extensions~\cite{MackaaySavage2018,Brundan2018} of Khovanov's Heisenberg category.
These generalizations have been incorporated by~\cite{BrundanSavageWebster2020} into the category
denoted by $\cH eis_{k}(z,t)$, where $z$ and $t$ are invertible elements of the ground field and $k\in \bZ$.
It is interesting to study if these quantum Heisenberg categories can give useful tools to study the asymptotic representation of some quantum deformed algebraic structures.

One candidate could be a family of cyclotomic Hecke algebras.
Given a polynomial
\begin{equation*}
    f(x)=f_{0}+f_{1}x+\cdots+f_{l}x^{l}
\end{equation*}
of degree $l$ such that $f_{0}f_{1}=1$, the cyclotomic Hecke algebras $H_{n}^{f}$, $n\in \bZ_{\geq 0}$ are defined as quotients of affine Hecke algebras.
As an analogue of the functors $\cF\colon \cH\to \cE nd(\cC)$ that we explained in Subsection~\ref{subsect:rep_of_H}, there is a functor~\cite{BrundanSavageWebster2020}
\begin{equation*}
    \cH eis_{-l}(z,f_{0}^{-1})\to \cE nd \left(\bigoplus_{n\geq 0}H^{f}_{n}-\catMod\right)
\end{equation*}
under the matching $z=q-q^{-1}$ with $q$ being the standard quantization parameter of Hecke algebras.
Then it is an interesting problem to ask if the construction of \cite{KLM2019} can be extended to this setting.

Another direction of deformation is the Jack deformation.
Although the graphical realizations of Jack polynomials are still absent,
the Jack inner product has been graphically realized by~\cite{LicataRossoSavage2018} in the context of Frobenius Heisenberg categorification~\cite{CautisLicata2012,RossoSavage2017} (see also~\cite{Savage2019} for the Frobenius Heisenberg categorification).
One can ask if the results therein are useful in the context of asymptotic representation theory.
It is especially natural to expect that the Jack deformed normalized characters are polynomials of suitably defined Boolean cumulants with non-negative integer coefficients.
In this direction, it is also important to look for the graphical realization of free cumulants.
In fact, there are several attempts to understand the Jack deformed analogue of Kerov polynomials~\cite{Lassalle2009,DolegaFeray2016}.

The analogue of \cite{KLM2019} associated with projective (spin) representations of symmetric groups has been established in~\cite{KvingeOguzReeks2020} using the twisted Heisenberg category~\cite{CautisSussan2015,OguzReeks2017}.
It is plausible that our main results have analogues in the case of projective representations.
We will report some results along this line in the forthcoming paper.
Note that spin normalized characters and the spin analogue of Kerov polynomials have been also studied in~\cite{Matsumoto2018,MatsumotoSniady2020}.

\bibliographystyle{alpha}
\bibliography{heis_cat}

\end{document}